\documentclass[leqno,12pt]{article}

\usepackage[dvips]{graphics}
\usepackage{color}
\usepackage{ifthen}
\usepackage{soul}
\usepackage{latexsym}
\usepackage{amsthm,amssymb}
\usepackage{amsbsy,amsfonts,amsmath}

\newcommand{\bysame}{\leavevmode\hbox to 3em{\hrulefill}\,}

\usepackage{a4}

\numberwithin{equation}{section}

\newtheorem{theorem}{Theorem}[section]
\newtheorem{proposition}[theorem]{Proposition}

\newtheorem{lemma}[theorem]{Lemma}
\newtheorem{example}[theorem]{Example}

\theoremstyle{definition}

\theoremstyle{remark}
\newtheorem{remark}[theorem]{Remark}

\newcommand{\al}{\alpha}

\newcommand{\bC}{{\mathbb{C}}}

\newcommand{\krdel}{\delta}

\newcommand{\pinfty}{\infty}

\newcommand{\bR}{\mathbb{R}}
\newcommand{\fkJ}{J}
\newcommand{\bZ}{\mathbb{Z}}
\newcommand{\bZgeqo}{\bZ_{\geq 0}}
\newcommand{\bRgeqo}{\bR_{\geq 0}}
\newcommand{\bRlo}{\bR_{>0}}

\newcommand{\bN}{\mathbb{N}}

\newcommand{\bK}{\mathbb{K}}
\newcommand{\bKt}{\bK^\times}

\newcommand{\kpch}{{\hat{o}}}

\newcommand{\fkA}{\mathfrak{A}}
\newcommand{\bhm}{\chi}

\newcommand{\rmrank}{{\mathrm{rank}}}

\newcommand{\tfkA}{\fkA}
\newcommand{\tPi}{\Pi}

\newcommand{\tfkAPi}{\tfkA}

\newcommand{\fkI}{I}

\newcommand{\tbhm}{\bhm}

\newcommand{\tlambda}{\lambda}
\newcommand{\tmu}{\mu}

\newcommand{\rmEnd}{{\mathrm{End}}}

\newcommand{\tL}{\trL} 

\newcommand{\rmSpan}{{\mathrm{Span}}}

\newcommand{\tX}{{\tilde X}}

\newcommand{\rmid}{{\rm {id}}}

\newcommand{\trUp}{U^+}
\newcommand{\trUm}{U^-}

\newcommand{\trK}{K}
\newcommand{\trL}{L}
\newcommand{\trE}{E}
\newcommand{\trF}{F}

\newcommand{\trDelta}{\Delta}
\newcommand{\trS}{S}
\newcommand{\trvarepsilon}{\varepsilon}

\newcommand{\prtrZnewparatbhmtPi}{\prtrZnewparatbhmtPi(\tbhm,\tPi)}

\newcommand{\trvt}{\vartheta}

\newcommand{\rmIm}{{\mathrm{Im}}}

\newcommand{\mclL}{{\mathcal{L}}}

\newcommand{\DPA}{{\mathsf{A}}}
\newcommand{\DPAcop}{\Delta_\DPA}
\newcommand{\DPAant}{S_\DPA}
\newcommand{\DPAcu}{\varepsilon_\DPA}

\newcommand{\DPB}{{\mathsf{B}}}
\newcommand{\DPBcop}{\Delta_\DPB}
\newcommand{\DPBant}{S_\DPB}
\newcommand{\DPBcu}{\varepsilon_\DPB}

\newcommand{\Dpli}{{\mathsf{p}}}

\newcommand{\DPBop}{\DPB^{\mathrm{op}}}

\newcommand{\Dpbi}{{\mathsf{P}}}

\newcommand{\DPD}{{\mathsf{D}}}
\newcommand{\DPDcop}{\Delta_\DPD}
\newcommand{\DPDant}{S_\DPD}
\newcommand{\DPDcu}{\varepsilon_\DPD}

\newcommand{\lact}{\bullet_{\mathrm{l}}}
\newcommand{\ract}{\bullet_{\mathrm{r}}}

\newcommand{\duDPD}{\DPD^*}

\newcommand{\DPDOm}{\Omega_\DPD}

\newcommand{\duact}{\cdot}

\newcommand{\brvfkI}{{\breve{\fkI}}}

\newcommand{\brvtrvt}{{\breve{\trvt}}}

\newcommand{\brva}{{\breve{a}}}
\newcommand{\brvVlm}{{\breve{V}}}
\newcommand{\brvvlm}{{\breve{v}}}

\newcommand{\brvm}{{\breve{m}}}

\newcommand{\brvrho}{{\breve{\rho}}}

\newcommand{\rmTrbrvVlm}{{{\mathrm{Tr}}_\brvVlm}}

\newcommand{\brvf}{{\breve{f}}}

\newcommand{\brvY}{{\breve{Y}}}

\newcommand{\brvt}{{\breve{t}}}

\newcommand{\brvVlmdual}{{\brvVlm^*}}

\newcommand{\BYcnal}{\alpha}
\newcommand{\BYcnl}{l}

\newcommand{\BYcnfkJ}{J}

\newcommand{\BYcnkpch}{{\hat{o}}}

\newcommand{\BYcnfkAoriginal}{{\mathfrak{A}}}

\newcommand{\BYcnchi}{\chi}
\newcommand{\BYcnfkI}{I}
\newcommand{\BYcnpi}{\pi}
\newcommand{\BYcnpial}{\al}
\newcommand{\BYcnfkAoriginalpi}{\BYcnfkAoriginal_\BYcnpi}
\newcommand{\BYcnfkAoriginalpip}{\BYcnfkAoriginalpi^+}
\newcommand{\BYcnU}{U}
\newcommand{\BYcnUchipi}{\BYcnU(\BYcnchi,\BYcnpi)}
\newcommand{\BYcntrK}{K}
\newcommand{\BYcntrL}{L}
\newcommand{\BYcntrE}{E}
\newcommand{\BYcntrF}{F}

\newcommand{\BYcnUdual}{\BYcnU^*}

\newcommand{\BYcnrmSpan}{{\mathrm{Span}}}

\newcommand{\BYcnUo}{\BYcnU^0}
\newcommand{\BYcnUochipi}{\BYcnUo(\BYcnchi,\BYcnpi)}
\newcommand{\BYcnUp}{\BYcnU^+}
\newcommand{\BYcnUpchipi}{\BYcnUp(\BYcnchi,\BYcnpi)}
\newcommand{\BYcnUm}{\BYcnU^-}
\newcommand{\BYcnUmchipi}{\BYcnUm(\BYcnchi,\BYcnpi)}

\newcommand{\BYcntrm}{\varsigma}
\newcommand{\BYcntrmone}{\BYcntrm_1}
\newcommand{\BYcntrmtwo}{\BYcntrm_2}

\newcommand{\BYcnDelta}{\Delta}
\newcommand{\BYcnS}{S}
\newcommand{\BYcne}{\varepsilon}
\newcommand{\BYcnvartheta}{\vartheta}
\newcommand{\BYcnvarthetachipi}{\BYcnvartheta^{\BYcnchi,\BYcnpi}}

\newcommand{\BYcnUpflat}{\BYcnU^{+,\flat}}
\newcommand{\BYcnUmflat}{\BYcnU^{-,\flat}}

\newcommand{\BYcnUopflat}{\BYcnU^{+,\flat,0}}
\newcommand{\BYcnUomflat}{\BYcnU^{-,\flat,0}}

\newcommand{\BYcnrmid}{{\mathrm{id}}}

\newcommand{\DPX}{{\mathsf{X}}}

\newcommand{\BYcnSh}{{\mathfrak{Sh}}}
\newcommand{\BYcnShchipi}{\BYcnSh^{\BYcnchi,\BYcnpi}}

\newcommand{\BYcnbrvq}{\xi}

\newcommand{\BYcnRchipip}{{R^{\BYcnpi,+}_\BYcnchi}}
\newcommand{\BYcnvphpip}{{\varphi^{\BYcnpi,+}_\BYcnchi}}

\newcommand{\BYcnomega}{\omega}

\newcommand{\BYcnprtrZ}{\mathfrak{Z}}
\newcommand{\BYcnprtrZomchipi}{\BYcnprtrZ_\BYcnomega(\BYcnchi,\BYcnpi)}

\newcommand{\BYcntrRtbhmtPi}{{{\mathfrak{B}}^{\BYcnchi,\BYcnpi}_\BYcnomega}}

\newcommand{\BYcnq}{q}
\newcommand{\BYcnqbeta}{\BYcnq_\beta}
\newcommand{\BYcncbeta}{c_\beta}

\newcommand{\BYcnomegachhilmb}{{\BYcnomega^\BYcnchi_{\lambda,\mu;\beta}}}

\newcommand{\BYcnhrhochipi}{{\hat{\rho}}^{\BYcnchi,\BYcnpi}}
\newcommand{\BYcnhrhochipibeta}{{\BYcnhrhochipi(\beta)}}
\newcommand{\BYcnHCbhmpiomega}{{\mathfrak{HC}}^{\BYcnchi,\BYcnpi}_\BYcnomega}

\newcommand{\BYcnomegabrzero}{{\breve{\omega}}_0}
\newcommand{\BYcnprtrZomchipibrzero}{\BYcnprtrZ_{\BYcnomegabrzero}(\BYcnchi,\BYcnpi)}

\newcommand{\BYcnHCbhmpiomegabrzero}{{\mathfrak{HC}}^{\BYcnchi,\BYcnpi}_{\BYcnomegabrzero}}

\newcommand{\BYcnDPDBYcnU}{\Omega_\BYcnU} 

\newcommand{\BYcnUpi}{\BYcnU_\BYcnpi}
\newcommand{\BYcnUpichipi}{\BYcnUpi(\BYcnchi,\BYcnpi)}
\newcommand{\BYcnUopi}{\BYcnUo_\BYcnpi}
\newcommand{\BYcnUopichipi}{\BYcnUopi(\BYcnchi,\BYcnpi)}

\newcommand{\BYcnLamlmnewpr}{\Lambda^\BYcnchi_{\lambda,\mu;\BYcnomega}}

\newcommand{\BYcnFinnewpr}{{\mathrm{Fin}}^\BYcnchi_\BYcnomega}

\newcommand{\BYcnFinnewprZlm}{Z^\BYcnchi_{\lambda,\mu;\BYcnomega}}

\newcommand{\sangaq}{q}
\newcommand{\sangazeta}{\zeta}

\newcommand{\sangaV}{V}

\newcommand{\bXI}{{{\mathbb{X}}(I)}}
\newcommand{\bXsimI}{{{\mathbb{X}}^\sim(I)}}
\newcommand{\bLSsimI}{{{\mathbb{LS}}^\sim(I)}}

\newcommand{\acuteBYcnUpflat}{{\acute{U}}^{+,\flat}}
\newcommand{\acuteBYcntrK}{{\acute K}}
\newcommand{\acuteBYcntrE}{{\acute E}}

\newcommand{\acuteBYcnUpflato}{{\acute{U}}^{+,\flat,0}}
\newcommand{\acuteBYcnUp}{{\acute{U}}^+}
\newcommand{\acuteG}{{\acute{G}}}
\newcommand{\bB}{{\mathbb{B}}}

\newcommand{\acuteBYcntrm}{{\acute{\varsigma}}}

\newcommand{\bH}{{\mathbb{H}}}

\newcommand{\acutebH}{{\acute{\bH}}}
\newcommand{\acuteh}{{\acute{h}}}

\newcommand{\acuteRchip}{{{\acute{R}}_\chi^+}}
\newcommand{\acutephichip}{{{\acute{\varphi}}_\chi^+}}

\begin{document}


\title{Natural Elements of \\ Center of Generalized Quantum Groups}
\author{Punita Batra and Hiroyuki Yamane}
\date{}







\maketitle


\footnote[0]{{\tt{THIS VERSION IS ON 17/MAY/2019.}} \par
\par HY is partially supported by JSPS Grand-in-Aid for Scientific Research (C), 16K05095.}


\maketitle

\begin{abstract}
This paper gives elements in the (skew)  center of
the generalized quantum group corresponding to its irreducible finite dimensional modules.
Finally we give a conjecture stating that those must form a basis of the center.
\end{abstract}

\section{Introduction}
This paper is a continuation of \cite{BY15}, \cite{BY18}.
Let $V$ be a finite-dimensional real linear space.
Let $l:=\dim_\bR V$, and assume $l\geq 1$.
Let $\{\al_1,\ldots,\al_l\}$ be an $\bR$-basis of $V$.
So $V=\oplus_{i=1}^l\bR\al_i$.
Let $V_\bZ:=\oplus_{i=1}^l\bZ\al_i$.
Then $V_\bZ$ is a free $\bZ$-module with $\rmrank_\bZ V_\bZ=l$.
Let $\bK$ be a field. Let $\bKt:=\bK\setminus\{0\}$.
Let $\BYcnchi:V_\bZ\times V_\bZ\to\bKt$
be a map such that
$\BYcnchi(\lambda+\mu,\nu)=\BYcnchi(\lambda,\nu)\BYcnchi(\mu,\nu)$
and $\BYcnchi(\lambda,\mu+\nu)=\BYcnchi(\lambda,\mu)\BYcnchi(\lambda,\nu)$.
For the $\chi$, 
in the same way as the Lusztig way~\cite[3.1.1]{b-Lusztig93} to define the
quantum groups $U_q({\mathfrak{g}})$,
we can define the Hopf algebra $\BYcnU=\BYcnUchipi$ over $\bK$, 
which we call {\it{the generalized quantum group}}.
($\BYcnU$ is the one of  Subsection~\ref{subsection:DEFtU}
with assuming $V_\bZ$ to be $\BYcnfkAoriginal$. In Introduction, we may also assume  $V_\bZ=
\BYcnfkAoriginal=\BYcnfkAoriginalpi$. The notation $\pi$
means the map from $\{1,\ldots,l\}$ to $\BYcnfkAoriginal$
defined by $\pi(i):=\al_i$.)
For some $\BYcnchi$, 
$\BYcnU$ can be the quantum groups, the small quantum groups
at root of $1$, the quantum superalgebras, or the one associated with the Nichols algebras
in the Heckenberger's list \cite{Hec09}.
Further studies concerning $\BYcnU$ have been achieved by 
\cite{AY15}, \cite{AY18}, \cite{AYY12}, \cite{CH15}, \cite{HY08}, \cite{HY10}, \cite{JMY18},  etc.
Let $\BYcnU=\BYcnUm\otimes\BYcnUo\otimes\BYcnUp$ be the triangular
decomposition (see $(\BYcnU 5)$).
Let $\BYcnU=\oplus_{\lambda\in V_\bZ}\BYcnU_\lambda$ be the $V_\bZ$-grading (see $(\BYcnU 4)$).
Let $\BYcnRchipip (\subset V_\bZ)$ be the Kharchenko's positive root system defined for $\BYcnU$
(see Theorem~\ref{theorem:KhPBW}, which was originally given by \cite{Kha99}).
As in the assumption of Theorem~\ref{theorem:previousmain},
assume that $\BYcnRchipip$
is a finite set and assume that $\BYcnchi(\al,\al)\ne 1$ for all $\al\in\BYcnRchipip$.
Let $\BYcnomega:V_\bZ\to\bKt$ be the map such that
$\BYcnomega(\lambda+\mu)=\BYcnomega(\lambda)\BYcnomega(\mu)$.
Let
$\BYcnprtrZomchipi:=\{\,Z\in\BYcnU_0\,|\,\forall \lambda\in V_\bZ,\,\forall X\in\BYcnU_\lambda,\,
ZX
=\BYcnomega(\lambda)XZ\,\}$ (see \eqref{eqn:defprtrZomega}).
In \cite[Theorem~10.4]{BY18}  (see Theorem~\ref{theorem:previousmain}),
we have the Harish-Chandra type isomorphism
$\BYcnHCbhmpiomega:\BYcnprtrZomchipi\to\BYcntrRtbhmtPi$,
where $\BYcntrRtbhmtPi$ is the $\bK$-subspace of $\BYcnUo$
defined by the equrations
$(e1)_\beta$-$(e4)_\beta$ for all $\beta\in\BYcnRchipip$. 
As in \eqref{eqn:mainweight}, for $\lambda$, $\mu\in V_\bZ$,
define the $\bK$-algebra homomorphism
$\BYcnLamlmnewpr:\BYcnUo\to\bK$
by
$\BYcnLamlmnewpr(\trK_{\lambda^\prime}\trL_{\mu^\prime}):
=\BYcnchi(\lambda,\mu^\prime)
\BYcnchi(\lambda^\prime,\mu)\BYcnomega(\lambda^\prime)$
for $\lambda^\prime$, $\mu^\prime\in V_\bZ$,
where $\{K_{\lambda^\prime}L_{\mu^\prime}|\lambda^\prime, \mu^\prime\in V_\bZ\}$
is a $\bK$-basis of $\BYcnUo$ (see $(\BYcnU 3)$).
As in \eqref{eqn:mainweightd}, let
$\BYcnFinnewpr:=\{(\lambda,\tmu)\in V_\bZ\times V_\bZ|\dim\mclL(\BYcnLamlmnewpr)<\infty\}$,
where $\mclL(\BYcnLamlmnewpr)$ is the finite-dimensional simple
$\BYcnU$-modules defined in the way that $\BYcnLamlmnewpr$
is regarded as its highest weight.
By our main result Theorem~\ref{theorem:mainth}, for each $(\lambda,\tmu)\in\BYcnFinnewpr$,
we have $\BYcnFinnewprZlm\in\BYcnprtrZomchipi$
so that $\BYcnHCbhmpiomega(\BYcnFinnewprZlm)$ can be viewed as the character of
$\mclL(\BYcnLamlmnewpr)$.
As a final stage of this paper, we state Conjecture~\eqref{eqn:Conjecture}, which states that 
$\{\,\BYcnFinnewprZlm\,|\,(\lambda,\mu)\in\BYcnFinnewpr\,\}$
is a $\bK$-basis of $\BYcnprtrZomchipi$.

\section{Generalized quantum groups}\label{section:GQG}

\subsection{Preliminary}\label{subsection:Prel}

For $x$, $y\in\bR$,
let $\BYcnfkJ_{x,y}:=\{z\in\bZ|x\leq z\leq y\}$.
Let $\bK$ be a field. Let $\bKt:=\bK\setminus\{0\}$.
For $n\in\bZgeqo$ and $x\in\bK$, let $(n)_x:=\sum_{r=1}^n x^{r-1}$,
and $(n)_x!:=\prod_{r=1}^n(r)_x$.
For $n\in\bZgeqo$, $m\in\fkJ_{0,n}$ and $x\in\bK$,
define ${n\choose m}_x\in\bK$ by
${n\choose 0}_x:={n\choose n}_x:=1$,
and ${n\choose m}_x:={n-1\choose m}_x+x^{n-m}{n-1\choose m-1}_x
=x^m{n-1\choose m}_x+{n-1\choose m-1}_x$
(if $m\in\fkJ_{1,n-1}$).
If $(m)_x!(n-m)_x!\ne 0$, then ${n\choose m}_x=
{\frac {(n)_x!} {(m)_x!(n-m)_x!}}$.
For $x$, $y$, $z\in\bK$, and $n\in\bN$, we have
$\prod_{t=0}^{n-1}(y+x^tz)=\sum_{m=0}^n
x^{{\frac {m(m-1)} 2}}{n\choose m}_xy^{n-m}z^m$.

For $x\in\bKt$, define $\BYcnkpch(x)\in\bZgeqo\setminus\{1\}$ by
\begin{equation}\label{eqn:chd}
\BYcnkpch(x):= 
\left\{\begin{array}{l}\min\{\,r^\prime\in\fkJ_{2,\pinfty}\,|\,(r^\prime)_x!=0\,\} 
\,\,\mbox{if $(r^{\prime\prime})_x!=0$ for some $r^{\prime\prime}\in\fkJ_{2,\pinfty}$}, \\
0 \quad\mbox{otherwise.}
\end{array}\right.
\end{equation}

Let $\BYcnfkAoriginal$ be an abelian group. Then $\BYcnfkAoriginal$
is a $\bZ$-module.
Let $\BYcnchi:\BYcnfkAoriginal\times\BYcnfkAoriginal\to\bKt$ be a map such that
\begin{equation}\label{eqn:bich}
\forall \lambda, \forall \mu, \forall \nu\in\BYcnfkAoriginal,\,\BYcnchi(\lambda+\mu,\nu)=\BYcnchi(\lambda,\nu)\BYcnchi(\mu,\nu),\,
\BYcnchi(\lambda,\mu+\nu)=\BYcnchi(\lambda,\mu)\BYcnchi(\lambda,\nu).
\end{equation}

Let $\BYcnl\in\bN$ and $\BYcnfkI:=\BYcnfkJ_{1,\BYcnl}$.
Assume that there exists  an injection $\BYcnpi:\BYcnfkI\to\BYcnfkAoriginal$ 
such that $\BYcnpi(\BYcnfkI)$ is a $\bZ$-basis of $\BYcnrmSpan_\bZ(\BYcnpi(\BYcnfkI))$.
Let $\BYcnfkAoriginalpi:=\BYcnrmSpan_\bZ(\BYcnpi(\BYcnfkI))$.
Let $\BYcnpial_i:=\BYcnpi(i)$ ($i\in\BYcnfkI$).
Then $\BYcnfkAoriginalpi:=\oplus_{i\in\BYcnfkI}\bZ\BYcnpial_i$.
Let $\BYcnfkAoriginalpip:=\oplus_{i\in\BYcnfkI}\bZgeqo\BYcnpial_i$.

\subsection{Definition of $\BYcnU=\BYcnUchipi$}\label{subsection:DEFtU}

The facts mentioned in this subsection is well-known
and can be proved in a standard way introduced by Drinfeld~\cite{Dr86}
and Lusztig~\cite[3.1.1]{b-Lusztig93}.

There exists a unique associative $\bK$-algebra (with $1$)
$\BYcnU=\BYcnUchipi$ satisfying the following conditions $(\BYcnU 1)$-$(\BYcnU 6)$.
\newline\newline
$(\BYcnU 1)$ As a $\bK$-algebra, $\BYcnU$ is generated by the elements:
\begin{equation}\label{eqn:gene}
\BYcntrK_\lambda,\,\BYcntrL_\lambda\,\,(\lambda\in\BYcnfkAoriginal),\quad
\BYcntrE_i, \BYcntrF_i\,\,(i\in\BYcnfkI).
\end{equation} \newline
$(\BYcnU 2)$ The elements of \eqref{eqn:gene} satisfy the following relations.
\begin{equation}\label{eqn:relone}
\begin{array}{l}
\BYcntrK_0=\BYcntrL_0=1,\,
\BYcntrK_\lambda \BYcntrK_\mu=\BYcntrK_{\lambda+\mu},\,
\BYcntrL_\lambda \BYcntrL_\mu=\BYcntrL_{\lambda+\mu},\,
\BYcntrK_\lambda \BYcntrL_\mu=\BYcntrL_\mu \BYcntrK_\lambda, \\
\BYcntrK_\lambda \BYcntrE_i =\BYcnchi(\lambda,\BYcnpial_i)\BYcntrE_i \BYcntrK_\lambda,\,
\BYcntrL_\lambda \BYcntrE_i  =\tbhm(-\BYcnpial_i,\lambda)\BYcntrE_i \BYcntrL_\lambda,\\
\BYcntrK_\lambda \BYcntrF_i  =\tbhm(\lambda,-\BYcnpial_i)\BYcntrF_i \BYcntrK_\lambda,\,
\BYcntrL_\tlambda \BYcntrF_i  =\tbhm(\BYcnpial_i,\lambda)\BYcntrF_i\BYcntrL_\lambda,\\
\mbox{$[\BYcntrE_i,\BYcntrF_j]=\delta_{ij}(-\BYcntrK_{\BYcnpial_i}+\BYcntrL_{\BYcnpial_i})$}.
\end{array}
\end{equation}\newline
$(\BYcnU 3)$ Define the map $\BYcntrmone:\BYcnfkAoriginal\times\BYcnfkAoriginal\to\BYcnU$
by $\BYcntrmone(\lambda,\mu):=\BYcntrK_\lambda\BYcntrL_\mu$.
Define the $\bK$-subalgebra $\BYcnUo=\BYcnUochipi$ of $\BYcnU$ by 
$\BYcnUo:=\BYcnrmSpan_\bK(\BYcntrmone(\BYcnfkAoriginal\times\BYcnfkAoriginal))$.
Then  $\BYcntrmone$ is injective,
and $\BYcntrmone(\BYcnfkAoriginal\times\BYcnfkAoriginal)$ is a $\bK$-basis of  $\BYcnUo$.
\newline\newline
$(\BYcnU 4)$ There exist $\bK$-subspaces $\BYcnU_\lambda=\BYcnUchipi_\lambda$ of $\BYcnU$
for $\lambda\in\BYcnfkAoriginalpi$ satisfying the following conditions $(\BYcnU 4-1)$-$(\BYcnU 4-3)$.
\newline\newline
$(\BYcnU 4-1)$ We have $\BYcnUo\subset\BYcnU_0$ and $\BYcntrE_i\in\BYcnU_{\BYcnpial_i}$,
$\BYcntrF_i\in\BYcnU_{-\BYcnpial_i}$ ($i\in\BYcnfkI$). \newline
$(\BYcnU 4-2)$ We have $\BYcnU_\lambda\BYcnU_\mu\subset
\BYcnU_{\lambda+\mu}$ ($\lambda$, $\mu\in\BYcnfkAoriginalpi$).  \newline
$(\BYcnU 4-3)$ We have 
$\BYcnU=\oplus_{\lambda\in\BYcnfkAoriginalpi}\BYcnU_\lambda$
as a $\bK$-linear spaces.
\newline\newline
$(\BYcnU 5)$ Let $\BYcnUp=\BYcnUpchipi$ 
(resp. $\BYcnUm=\BYcnUmchipi$) be the $\bK$-subalgebra (with $1$) of $\BYcnU$
generated by $\BYcntrE_i$ (resp. $\BYcntrF_i$) ($i\in\BYcnfkI$).
Define the $\bK$-linear homomorphism $\BYcntrmtwo:\BYcnUm\otimes_\bK\BYcnUo\otimes_\bK\BYcnUp\to\BYcnUchipi$
by $\BYcntrmtwo(Y\otimes Z\otimes X):=YZX$.
Then $\BYcntrmtwo$ is a $\bK$-linear isomorphism.
\newline\newline
$(\BYcnU 6)$  For $\lambda\in\BYcnfkAoriginalpi$,
define $\BYcnUp_\lambda=\BYcnUpchipi_\lambda$
(resp. $\BYcnUm_\lambda=\BYcnUmchipi_\lambda$)
by $\BYcnUp_\lambda:=\BYcnUp\cap\BYcnU_\lambda$
(resp. $\BYcnUm_\lambda=\BYcnUm\cap\BYcnU_\lambda$).
Then for $\lambda\in\BYcnfkAoriginalpip\setminus\{0\}$,
we have 
$\{X\in\BYcnUp_\lambda|\forall i\in\BYcnfkI, [X,\BYcntrF_i]=0\}=\{0\}$
and $\{Y\in\BYcnUm_{-\lambda}|\forall i\in\BYcnfkI, [\BYcntrE_i,Y]=0\}=\{0\}$.
\newline\newline
Notice that $\BYcnUp=\oplus_{\lambda\in\BYcnfkAoriginalpip}\BYcnUp_\lambda$,
$\BYcnUm=\oplus_{\lambda\in\BYcnfkAoriginalpip}\BYcnUm_{-\lambda}$,
$\BYcnUp_0=\BYcnUm_0=\bK\cdot1_\BYcnU$
and $\BYcnUp_{\BYcnal_i}=\bK\cdot\BYcntrE_i$,
$\BYcnUm_{-\BYcnpial_i}=\bK\cdot\BYcntrF_i$
($i\in\BYcnfkI$).
\newline\par
We also regard $\BYcnU=\BYcnUchipi$ as a Hopf algebra $(\BYcnU,\BYcnDelta,\BYcnS,\BYcne)$
by
\begin{equation}\label{eqn:defHopf}
\begin{array}{l}
\BYcnDelta(\BYcntrK_\lambda)=\BYcntrK_\lambda\otimes\BYcntrK_\lambda,
\BYcnDelta(\BYcntrL_\lambda)=\BYcntrL_\lambda\otimes\BYcntrL_\lambda,
\BYcnDelta(\BYcntrE_i)=\BYcntrE_i\otimes  1+\BYcntrK_{\BYcnpial_i}\otimes\BYcntrE_i, \\
\BYcnDelta(\BYcntrF_i)=\BYcntrF_i\otimes
\BYcntrL_{\BYcnpial_i}+1\otimes  \BYcntrF_i,
\BYcnS(\BYcntrK_\lambda)=\BYcntrK_{-\lambda}, 
\BYcnS(\BYcntrL_\lambda)=\BYcntrL_{-\lambda}, \\
\BYcnS(\BYcntrE_i)=-\BYcntrK_{-\BYcnpial_i}\BYcntrE_i,
\BYcnS(\BYcntrF_i)=-\BYcntrF_i\BYcntrL_{-\BYcnpial_i}, \\
\BYcne(\BYcntrK_\lambda)=\BYcne(\BYcntrL_\lambda)=1,
\BYcne(\BYcntrE_i)=\BYcne(\BYcntrF_i)=0.
\end{array}
\end{equation}

Let $\BYcnUpflat:=\oplus_{\lambda\in\BYcnfkAoriginal}\BYcnUp\BYcntrK_\lambda$,
and $\BYcnUmflat:=\oplus_{\lambda\in\BYcnfkAoriginal}\BYcnUm\BYcntrL_\lambda$.
Then $\BYcnU=\BYcnrmSpan_\bK(\BYcnUmflat\BYcnUpflat)=\BYcnrmSpan_\bK(\BYcnUpflat\BYcnUmflat)$.

As in a
standard way (see \cite{Dr86}), we have a bilinear form
$\BYcnvartheta=\BYcnvarthetachipi:\BYcnUpflat\times\BYcnUmflat\to\bK$ having the following
properties, see also Subsection~\ref{subsection:Hp} below.
\begin{equation}\label{eqn:bprtoftB}
\begin{array}{l}
\BYcnvartheta(\BYcntrK_\lambda,\BYcntrL_\mu)=\BYcnchi(\lambda,\mu),
\BYcnvartheta(\BYcntrE_i,\BYcntrF_j)=\delta_{ij},
\BYcnvartheta(\BYcntrK_\lambda,\BYcntrF_j)=\BYcnvartheta(\BYcntrE_i,\BYcntrL_\lambda)=0, \\
\BYcnvartheta(X^+Y^+,X^-)=
\sum_{k^-}\BYcnvartheta(X^+,(X^-)^{(2)}_{k^-})\BYcnvartheta(Y^+,(X^-)^{(1)}_{k^-}),\\
\BYcnvartheta(X^+,X^-Y^-)=
\sum_{k^+}\BYcnvartheta((X^+)^{(1)}_{k^+},X^-)\BYcnvartheta((X^+)^{(2)}_{k^+},Y^-),\\
\BYcnvartheta(\BYcnS(X^+),X^-)=\BYcnvartheta(X^+,\BYcnS^{-1}(\tX^-)),\\
\BYcnvartheta(X^+,1)=\BYcne(X^+),
\BYcnvartheta(1,X^-)=\BYcne(X^-), \\
X^-X^+  =\sum_{r^+,r^-}
\BYcnvartheta((X^+)^{\prime,(1)}_{r^+}, \BYcnS((X^-)^{\prime,(1)}_{r^-}))
\BYcnvartheta((X^+)^{\prime,(3)}_{r^+}, (X^-)^{\prime,(3)}_{r^-}) \\
\quad\quad\quad\quad\quad\quad\quad\quad
\cdot(X^+)^{\prime,(2)}_{r^+}(X^-)^{\prime,(2)}_{r^-}, \\
X^+X^- =\sum_{r^+,r^-}
\BYcnvartheta((X^+)^{\prime,(3)}_{r^+}, \BYcnS((X^-)^{\prime,(3)}_{r^-}))
\BYcnvartheta((X^+)^{\prime,(1)}_{r^+}, (X^-)^{\prime,(1)}_{r^-}) \\
\quad\quad\quad\quad\quad\quad\quad\quad
\cdot(X^-)^{\prime,(2)}_{r^-}(X^+)^{\prime,(2)}_{r^+}
\end{array}
\end{equation} for $\lambda$, $\mu\in\BYcnfkAoriginal$, $i$, $j\in\BYcnfkI$,
and $X^+$, $Y^+\in\BYcnUpflat$, $X^-$, $Y^-\in\BYcnUmflat$,
where $(X^+)^{(x)}_{k^+}$ and $(X^-)^{(x)}_{k^-}$ with $x\in\BYcnfkJ_{1,2}$
(resp. $(X^+)^{\prime,(y)}_{r^+}$ and
$(X^-)^{\prime,(y)}_{r^-}$ with $y\in\BYcnfkJ_{1,3}$)
are any elements of $\BYcnUpflat$ and $\BYcnUmflat$ respectively
satisfying $\BYcnDelta(X^\pm)=\sum_{k^\pm}(X^\pm)^{(1)}_{k^\pm}
\otimes (X^\pm)^{(2)}_{k^\pm}$,
(resp. $((\BYcnrmid_\BYcnU\otimes\BYcnDelta)\circ\BYcnDelta)(X^\pm)=
\sum_{r^\pm}(X^\pm)^{\prime,(1)}_{r^\pm}\otimes
(X^\pm)^{\prime,(2)}_{r^\pm}\otimes (X^\pm)^{\prime,(3)}_{r^\pm}$).
We have
\begin{equation} \label{eqn:protvt}
\BYcnvartheta(X^+\BYcntrK_\tlambda,X^-\tL_\tmu)=\tbhm(\tlambda,\tmu)\BYcnvartheta(X^+,X^-)
\quad (\lambda,\mu\in\BYcnfkAoriginal,\,X^+\BYcnUpflat,\,X^-\BYcnUmflat),
\end{equation} and
\begin{equation} \label{eqn:protvtb}
\BYcnvartheta(\BYcnUp_\lambda,\BYcnUm_{-\mu})=\{0\}\quad\mbox{if $\lambda\ne\mu$.}
\end{equation}

By \eqref{eqn:protvt}-\eqref{eqn:protvtb}, we can easily see:

\begin{lemma} \label{lemma:BYcniclmm}
Let $\BYcnchi:\BYcnfkAoriginal\times\BYcnfkAoriginal\to\bKt$
and $\BYcnpi:\BYcnfkI\to\BYcnfkAoriginal$ be as above.
Let $\BYcnfkAoriginal^\prime$ be a $\bZ$-submodule of 
$\BYcnfkAoriginal$ such that $\BYcnpi(I)\subset\BYcnfkAoriginal^\prime$.
Let $\bK^\prime$ be a subfield of $\bK$ such that
$\BYcnchi(\BYcnfkAoriginal^\prime,\BYcnfkAoriginal^\prime)\subset(\bK^\prime)^\times$.
Let $\eta:\BYcnfkAoriginal^\prime\to\BYcnfkAoriginal$ be the inclusion map.
Define the map $\BYcnchi^\prime:
\BYcnfkAoriginal^\prime\times\BYcnfkAoriginal^\prime\to(\bK^\prime)^\times$ 
by $\eta\circ\BYcnchi^\prime=\BYcnchi_{|\BYcnfkAoriginal^\prime\times\BYcnfkAoriginal^\prime}$.
Define the map $\pi^\prime:\BYcnfkI\to\BYcnfkAoriginal^\prime$
by $\eta\circ\pi^\prime=\pi$.
Then there exists a $\bK$-algebra monomorphism
$f:U(\BYcnchi^\prime,\pi^\prime)\otimes_{\bK^\prime}\bK\to\BYcnUchipi$
such that $f(\BYcntrK_\lambda\BYcntrL_\mu)
=\BYcntrK_\lambda\BYcntrL_\mu$
$(\lambda,\mu\in\BYcnfkAoriginal^\prime)$ and
$f(\BYcntrE_i)=\BYcntrE_i$,
$f(\BYcntrF_i)=\BYcntrF_i$
$(i\in\BYcnfkI)$.
\end{lemma}

Define the $\bK$-linear map 
$\BYcnShchipi:\BYcnUchipi\to\BYcnUochipi$
by 
\begin{equation*}
\BYcnShchipi(Y\BYcntrK_\lambda\BYcntrL_\mu X)
=\BYcne(Y)\BYcne(X)\BYcntrK_\lambda
\BYcntrL_\mu
\end{equation*} 
$(X\in\BYcnUp,\,Y\in\BYcnUm,\,\lambda,\,\mu\in\BYcnfkAoriginal)$.


\subsection{Kharchenko's PBW theorem}

\begin{theorem}{\rm{(}}Kharchenko's PBW theorem~{\rm{\cite[Theorem~2]{Kha99}}}, {\rm{\cite[Theorem~2.2]{Kh15}}},
see also {\rm{\cite[Theorem~3.14]{HLec08}}} and
Section~{\rm{\ref{section:Appendix}}}.{\rm{)}} \label{theorem:KhPBW} 
Keep the notation as above. Then there exists a unique pair of $(\BYcnRchipip,\BYcnvphpip)$
of a subset $\BYcnRchipip$ of $\BYcnfkAoriginalpip\setminus\{0\}$ and a map $\BYcnvphpip:\BYcnRchipip\to\bN$ satisfying the following. 
Let $X:=\{(\al,t)\in\BYcnRchipip\times\bN|t\in\BYcnfkJ_{1,\BYcnvphpip(\al)}\}$.
Define the map $z:X\to\BYcnRchipip$ by $z(\al,t):=\al$.
Let $Y$ be the set of maps $y:X\to\bZgeqo$
such that $|\{x\in X|y(x)\geq 1\}|<\infty$ and $(y(x))_{\BYcnchi(z(x),z(x))}!\ne 0$
for all $x\in X$. Then
\begin{equation*}
\forall\lambda\in\BYcnfkAoriginalpip,\,\dim\BYcnUpchipi_\lambda=|\{y\in Y|\sum_{x\in X}y(x)z(x)=\lambda\}|.
\end{equation*}
\end{theorem}

\begin{theorem}{\rm{(\cite[Proposition~1]{Hec06}, \cite[Theorem~4.9]{HY10}}}) 
If $|\BYcnRchipip|<\infty$, then $\BYcnvphpip(\BYcnRchipip)=\{1\}$.
\end{theorem}

\subsection{Skew centers}
Let $\BYcnomega:\BYcnfkAoriginalpi\to\bKt$ be a $\bZ$-module homomorphism.

\begin{equation}\label{eqn:defprtrZomega}
\BYcnprtrZomchipi:=\{\,Z\in\BYcnUchipi_0\,|\,\forall \lambda\in\BYcnfkAoriginalpi,\,\forall X\in\BYcnUchipi_\lambda,\,
ZX
=\BYcnomega(\lambda)XZ\,\}.
\end{equation}

Define the $\bZ$-module homomorphism $\BYcnhrhochipi:\BYcnfkAoriginalpi\to\bKt$ by
\begin{equation*}
\BYcnhrhochipi(\BYcnpial_j):=\BYcnchi(\BYcnpial_j,\BYcnpial_j)\quad(j\in\BYcnfkI),
\end{equation*} where $\BYcnpial_j:=\BYcnpi(j)$, as above.

For each $\beta\in\BYcnRchipip$, let $\BYcntrRtbhmtPi(\beta)$ be the $\bK$-linear
subspace of $\BYcnUochipi$ formed by
the elements 
\begin{equation*}
\sum_{(\lambda,\mu)\in\BYcnfkAoriginal^2}a_{(\lambda,\mu)}
\BYcntrK_\lambda\BYcntrL_\mu
\end{equation*} with $a_{(\lambda,\mu)}\in\bK$
satisfying the following equations $(e1)_\beta$-$(e4)_\beta$.
In $(e1)_\beta$-$(e4)_\beta$, let $\BYcnqbeta:=\BYcnchi(\beta,\beta)$,
$\BYcncbeta:=\BYcnkpch(\BYcnqbeta)$ and 
$\BYcnomegachhilmb
:=\BYcnomega(\beta)\cdot{\frac {\BYcnchi(\beta,\mu)} {\BYcnchi(\lambda,\beta)}}$.
\newline\par
$(e1)_\beta$ For $(\lambda,\mu)\in\BYcnfkAoriginal^2$ and $t\in\bZ\setminus\{0\}$,
if $\BYcnqbeta\ne 1$,
$\BYcncbeta=0$ 
and $\BYcnomegachhilmb=\BYcnqbeta^t$,
then the equation
$a_{(\lambda+t\beta,\mu-t\beta)}
=\BYcnhrhochipibeta^t\cdot a_{(\lambda,\mu)}$ holds.
\newline\par
$(e2)_\beta$ For $(\lambda,\mu)\in\BYcnfkAoriginal^2$, if $\BYcncbeta=0$ and 
$\BYcnomegachhilmb\ne\BYcnqbeta^t$ for all $t\in\bZ$, 
then the equation $a_{(\lambda,\mu)}=0$ holds.
\newline\par
$(e3)_\beta$ For $(\lambda,\mu)\in\BYcnfkAoriginal^2$,
if $\BYcnqbeta\ne 1$, $\BYcncbeta\geq 2$ and 
$\BYcnomegachhilmb=\BYcnqbeta^t$ for some $t\in\BYcnfkJ_{1,\BYcncbeta-1}$, 
the equation
\begin{equation*}
\begin{array}{l}
\displaystyle{\sum_{x=-\infty}^\infty}a_{(\lambda+(\BYcncbeta x+t)\beta,\mu-(\BYcncbeta x+t)\beta)}
\BYcnhrhochipibeta^{-(\BYcncbeta x+t)} \\
\quad =\displaystyle{\sum_{y=-\infty}^\infty}a_{(\lambda+\BYcncbeta y\beta,\mu-\BYcncbeta y\beta)}
\BYcnhrhochipibeta^{-\BYcncbeta y}
\end{array}
\end{equation*} holds.
\newline\par
$(e4)_\beta$ For $(\lambda,\mu)\in\BYcnfkAoriginal^2$, 
if  $\BYcncbeta\geq 2$ and 
$\BYcnomegachhilmb\ne\BYcnqbeta^m$ for all $m\in\bZ$, then the $\BYcncbeta-1$ equations
\begin{equation*}
\begin{array}{l}
\displaystyle{\sum_{x=-\infty}^\infty}a_{(\lambda+(\BYcncbeta x+t)\beta,\mu-(\BYcncbeta x+t)\beta)}
\BYcnhrhochipibeta^{-(\BYcncbeta x+t)} \\
\quad =\displaystyle{\sum_{y=-\infty}^\infty}a_{(\lambda+\BYcncbeta y\beta,\mu-\BYcncbeta y\beta)}
\BYcnhrhochipibeta^{-\BYcncbeta y} \\
(t\in\BYcnfkJ_{1,\BYcncbeta-1})
\end{array}
\end{equation*} hold.
Let
\begin{equation*}
\BYcntrRtbhmtPi:=\bigcap_{\beta\in\BYcnRchipip}\BYcntrRtbhmtPi(\beta).
\end{equation*}
\begin{theorem}\label{theorem:previousmain}
{\rm{(}}{\rm{\cite[Theorem~10.4]{BY18}}}{\rm{)}}
Assume $\BYcnfkAoriginal=\BYcnfkAoriginalpi$.
Assume $|\BYcnRchipip|<\infty$.
Assume that $\BYcnchi(\al,\al)\ne 1$ for all $\al\in\BYcnRchipip$.
Then we have the $\bK$-linear isomorphism $\BYcnHCbhmpiomega:\BYcnprtrZomchipi\to\BYcntrRtbhmtPi$ defined by $\BYcnHCbhmpiomega(X):=\BYcnShchipi(X)$.
\end{theorem}
The statement of \cite[Theorem~10.4]{BY18} has claimed that the above theorem holds if $\bK$ is an algebraically closed field.
However, by the argument in \cite{BY18},
we can easily see that it really holds for any field.
  
\section{Elements of $\BYcnprtrZomchipi$
via finite dimensional representations} \label{section:CSCent}
In Subsection~\ref{section:CSCent}, 
we use argument similar to that in 
\cite[Sections~2 and 3]{C3Tanisaki}.
\subsection{Hopf pairing} \label{subsection:Hp}
In Subsection~\ref{subsection:Hp}, 
let 
\begin{equation*}
\mbox{$\DPA=(\DPA,\DPAcop,\DPAant,\DPAcu)$
and $\DPB=(\DPB,\DPBcop,\DPBant,\DPBcu)$}
\end{equation*}
be Hopf algebras over $\bK$ with $\DPAant$ and $\DPBant$
being bijective, and assume that 
there exists a Hopf pairing
$\Dpli:\DPA\times\DPB\to\bK$,
that is, $\Dpli(a_1a_2,b)=\Dpli(a_1\otimes a_2,\DPBcop(b))$,
$\Dpli(a,b_1b_2)=\Dpli(\DPAcop(a),b_1\otimes b_2)$,
$\Dpli(\DPAant(a),b)=\Dpli(a,\DPBant(b))$,
$\Dpli(a,1_\DPB)=\DPAcu(a)$,
and $\Dpli(1_\DPA,b)=\DPBcu(b)$
hold for all $a$, $a_1$, $a_2\in\DPA$ and
all $b$, $b_1$, $b_2\in\DPB$,
where $\Dpli(a_1\otimes a_2,b_1\otimes b_2)
:=p(a_1,b_1)p(a_2,b_2)$.
Note that from $\DPB$, we obtain
the Hopf algebra
$\DPBop=(\DPB,y\circ\DPBcop,\DPBant^{-1},\DPBcu)$,
where we define the $\bK$-linear isomorphism
$y:\DPB\otimes\DPB\to\DPB\otimes\DPB$ by
$y(b_1\otimes b_2):=b_2\otimes b_1$,
i.e., $\DPBop$ is the {\it{opposite Hopf algebra}} of $\DPB$.

It is well-known that we have the Hopf algebra $\DPD=\DPD(\Dpli,\DPA,\DPB)
=(\DPD,\DPDcop,\DPDant,\DPDcu)$ 
satisfying the conditions below.
\newline\par
$({\rm{D1}})$ $\DPD$ include $\DPA$ and $\DPBop$
as Hopf subalgebras.
The $\bK$-linear homomorphism $\DPA\otimes\DPB\to\DPD$,
$a\otimes b\mapsto ab$, is bijective. \par
$({\rm{D2}})$ As a $\bK$-algebra, $\DPD$ has the multiplication such that
for $a\in\DPA$ and $b\in\DPB$,
$ba=\sum_{i,j}\Dpli(\DPAant^{-1}(a^{(1)}_i),b^{(3)}_j)\Dpli(a^{(3)}_i,b^{(1)}_j)a^{(2)}_ib^{(2)}_j$,
where $((\DPAcop\otimes\rmid)\circ\DPAcop)(a)=\sum_i
a^{(1)}_i\otimes a^{(2)}_i\otimes a^{(3)}_i$ and  
$((\DPBcop\otimes\rmid)\circ\DPBcop)(b)=\sum_j
b^{(1)}_j\otimes b^{(2)}_j\otimes b^{(3)}_j$; we also have
$ab=\sum_{i,j}\Dpli(\DPAant^{-1}(a^{(3)}_i),b^{(1)}_j)\Dpli(a^{(1)}_i,b^{(3)}_j)b^{(2)}_ja^{(2)}_i$.
\newline\newline
We can easily see that the $\bK$-linear homomorphism $\DPB\otimes\DPA\to\DPD$,
$b\otimes a\to ba$, is bijective.
Define the $\bK$-bilinear map 
$\Dpbi:\DPD\times\DPD\to\bK$ by 
\begin{equation*}
\Dpbi(a_1b_1,b_2a_2):=\Dpli(a_1,b_2)\cdot\Dpli(a_2,b_1)\quad (a_1,a_2\in\DPA,\,
b_1,b_2\in\DPB).
\end{equation*}
Define the left (resp. right) $\bK$-algebra action $\lact$ 
(resp. $\ract$) of $\DPD$ on $\DPD$ by 
\begin{equation*}
x\lact y := \sum_ix^{(1)}_iy\DPDant(x^{(2)}_i)
\quad ({\rm{resp.}}\,\,y\ract x := \sum_i\DPDant(x^{(1)}_i)yx^{(2)}_i),
\end{equation*} where $\DPDcop(x)=\sum_ix^{(1)}_i\otimes x^{(2)}_i$.
We have $x_1\lact (x_2\lact y)=x_1x_2\lact y$
(resp. $(y\ract x_1)\ract x_2=y\ract x_1x_2$)
for $x_1$, $x_2$, $y\in\DPD$.
We have
\begin{lemma} \label{lemma:lmDad} 
It follows that 
\begin{equation} \label{eqn:Dad}
\Dpbi(x\lact y_1,y_2)=\Dpbi(y_1,y_2\ract x)
\quad (x,y_1,y_2\in\DPD).
\end{equation}
\end{lemma}
{\it{Proof.}} (cf. \cite[Proposition~2.2.1]{C3Tanisaki})
Letting $\DPX$ be $\DPA$ or $\DPB$,
for $c\in \DPX$ and $r\in\bN$ with $r\geq 2$,
let 
\begin{equation*}
\begin{array}{l}
\sum_i c^{(1)}_i\otimes\cdots\otimes c^{(r-1)}_i\otimes c^{(r)}_i \\
\quad :=((\Delta_\DPX\otimes\underbrace{\rmid_\DPX\otimes\cdots\otimes\rmid_\DPX}_{r-2})
\circ\cdots\circ(\Delta_\DPX\otimes\rmid_\DPX)\circ\Delta_\DPX)(c).
\end{array}
\end{equation*}
(Note that $\DPDcop(a)=\sum_i a^{(1)}_i\otimes a^{(2)}_i$ ($a\in\DPA$),
and  $\DPDcop(b)=\sum_j b^{(2)}_i\otimes b^{(1)}_j$ ($b\in\DPB$).)

Let $a,{\acute a}, {\grave a}\in\DPA$ and  $b,{\acute b}, {\grave b}\in\DPB$. We have
\begin{equation*}
\begin{array}{l}
a\lact{\acute a}{\acute b}=\sum_{i,j}\Dpli(a^{(4)}_i,{\acute b}^{(3)}_j)\Dpli(\DPAant(a^{(2)}_i),{\acute b}^{(1)}_j)
a^{(1)}_i{\acute a}\DPAant(a^{(3)}_i){\acute b}^{(2)}_j, \\
{\grave b}{\grave a}\ract a=\sum_{i,j}\Dpli(a^{(1)}_i,{\grave b}^{(1)}_j)\Dpli(\DPAant(a^{(3)}_i),{\grave b}^{(3)}_j)
{\grave b}^{(2)}_j\DPAant(a^{(2)}_i){\grave a}a^{(4)}_i, \\
b\lact{\acute a}{\acute b}=\sum_{i,j}\Dpli({\acute a}^{(1)}_i,\DPBant^{-1}(b^{(4)}_j))\Dpli({\acute a}^{(3)}_i,b^{(2)}_j)
{\acute a}^{(2)}_ib^{(3)}_j{\acute b}\DPBant^{-1}(b^{(1)}_j), \\
{\grave b}{\grave a}\ract b=\sum_{i,j}\Dpli({\grave a}^{(3)}_i,\DPBant^{-1}(b^{(1)}_j))\Dpli({\grave a}^{(1)}_i,b^{(3)}_j)
\DPBant^{-1}(b^{(4)}_j){\grave b}b^{(2)}_j{\grave a}^{(2)}_i.
\end{array}
\end{equation*} Hence we have
\begin{equation*}
\begin{array}{l}
\Dpbi(a\lact{\acute a}{\acute b},{\grave b}{\grave a}) \\
\quad =\sum_{i,j}\Dpli(a^{(4)}_i,{\acute b}^{(3)}_j)\Dpli(\DPAant(a^{(3)}_i),{\acute b}^{(1)}_j)
\Dpli(a^{(1)}_i{\acute a}\DPAant(a^{(3)}_i),{\grave b})\Dpli({\grave a},{\acute b}^{(2)}_j) \\
\quad =\sum_i\Dpli(\DPAant(a^{(3)}_i){\grave a}a^{(4)}_i,{\acute b})\Dpli(a^{(1)}_i{\acute a}\DPAant(a^{(3)}_i),{\grave b}) \\
\quad =\sum_{i,j}\Dpli(\DPAant(a^{(3)}_i){\grave a}a^{(4)}_i,{\acute b})
\Dpli(a^{(1)}_i,{\grave b}^{(1)}_j)\Dpli(\DPAant(a^{(3)}_i),{\grave b}^{(3)}_j)
\Dpli({\acute a},{\grave b}^{(2)}_j)  \\
\quad = \Dpbi({\acute a}{\acute b},{\grave b}{\grave a}\ract a),
\end{array}
\end{equation*} and
\begin{equation*}
\begin{array}{l}
\Dpbi(b\lact{\acute a}{\acute b},{\grave b}{\grave a}) \\
\quad =\sum_{i,j}\Dpli({\acute a}^{(1)}_i,\DPBant^{-1}(b^{(4)}_j))\Dpli({\acute a}^{(3)}_i,b^{(2)}_j)
\Dpli({\acute a}^{(2)}_i,{\grave b})\Dpli({\grave a},b^{(3)}_j{\acute b}\DPBant^{-1}(b^{(1)}_j)) \\
\quad =\sum_j\Dpli({\acute a},\DPBant^{-1}(b^{(4)}_j){\grave b}b^{(2)}_j)
\Dpli({\grave a},b^{(3)}_j{\acute b}\DPBant^{-1}(b^{(1)}_j))
 \\
\quad =\sum_{i,j}\Dpli({\acute a},\DPBant^{-1}(b^{(4)}_j){\grave b}b^{(2)}_j)
\Dpli({\grave a}^{(1)}_i,b^{(3)}_j)\Dpli({\grave a}^{(2)}_i,{\acute b})
\Dpli({\grave a}^{(3)}_i,\DPBant^{-1}(b^{(1)}_j))
  \\
\quad = \Dpbi({\acute a}{\acute b},{\grave b}{\grave a}\ract b).
\end{array}
\end{equation*}
This completes the proof.
\hfill $\Box$
\newline\par
Let $\duDPD$ be the $\bK$-linear space formed by all $\bK$-linear homomorphisms
from $\DPD$ to $\bK$.
Define the right $\bK$-algebra action $\duact$ of $\DPD$ on $\duDPD$ by
\begin{equation*}
(f\duact x)(y):=f(x\lact y) \quad (f\in\duDPD, x,y\in\DPD).
\end{equation*}
Define the $\bK$-linear homomorphism 
$\DPDOm:\DPD\to\duDPD$ by
$(\DPDOm(x))(y):=\Dpbi(y,x)$
($x,y\in\DPD$). By \eqref{eqn:Dad}, we have
\begin{equation} \label{eqn:DadOm}
\DPDOm(y) \duact x=\DPDOm(y \ract x) \quad (x,y\in\DPD).
\end{equation}
\begin{lemma} \label{lemma:Omequ}
Assume that $\Dpli$ is non-degenerate.
Let $x$, $y\in\DPD$. Then
\begin{equation} \label{eqn:interOmequ}
\DPDcu(y) x=x\ract y\quad\Longleftrightarrow\quad\DPDOm(x) \duact y=\DPDcu(y)\DPDOm(x).
\end{equation}
\end{lemma}
{\it{Proof.}} Since  $\Dpli$ is non-degenerate, $\Dpbi$ is so.
Hence $\DPDOm$ is injective.
Then \eqref{eqn:interOmequ} follows from \eqref{eqn:DadOm}.
\hfill $\Box$ 

\subsection{Extension of $\BYcnU=\BYcnUchipi$} \label{subsection:ExtUtbhmtPi}
Keep the notation in Subsection~\ref{subsection:DEFtU}. Recall $\BYcnfkI=\BYcnfkJ_{1,\BYcnl}$
and $\BYcnpial_i=\BYcnpi(i)$ ($i\in\BYcnfkI$). 
Let $\brvfkI:=\fkJ_{1,2\BYcnl+2}$, i.e., $\brvfkI=\BYcnfkI\cup\BYcnfkJ_{\BYcnl+1,2\BYcnl+2}$.
Let $\BYcnpial_j\in\BYcnfkAoriginal$
($j\in\BYcnfkJ_{\BYcnl+1,2\BYcnl+2}$). 
Let $\BYcnomega:\BYcnfkAoriginalpi\to\bKt$ be a $\bZ$-module homomorphism.
\begin{equation} \label{eqn:asstilq}
\begin{array}{l}
\mbox{In Subsections~\ref{subsection:ExtUtbhmtPi} and \ref{subsection:cstSkewce},
we assume the following $({\rm{ass}}1)$-$({\rm{ass}}4)$.} \\ 
\mbox{$({\rm{ass}}1)$ There exists $\BYcnbrvq\in\bKt\setminus\{1\}$ with 
$\kpch(\BYcnbrvq)=0$.} \\
\mbox{$({\rm{ass}}2)$ $\BYcnfkAoriginal$ is a $\bZ$-module with
$\rmrank_\bZ\BYcnfkAoriginal=2\BYcnl+2$
and $\BYcnfkAoriginal=\oplus_{i\in\brvfkI}\bZ\BYcnpial_i$.}\\
\mbox{$({\rm{ass}}3)$ One has $\BYcnchi(\BYcnal_i,\BYcnal_{j+\BYcnl+1})
=\BYcnchi(\BYcnal_{j+\BYcnl+1},\BYcnal_i)=\BYcnbrvq^{\delta_{ij}}$} \\
\mbox{and $\BYcnchi(\BYcnal_{i+\BYcnl+1},\BYcnal_{j+\BYcnl+1})=1$
for $i$, $j\in\BYcnfkJ_{1,\BYcnl+1}$.} \\
\mbox{$({\rm{ass}}4)$ One has $\BYcnomega(\BYcnal_i)=\BYcnchi(\BYcnal_i,\BYcnal_{l+1})$
and $\BYcnchi(\BYcnal_{l+1},\BYcnal_i)=1$ for $i\in\BYcnfkI$.}

\end{array}
\end{equation}

Let $X$ and $Y$ be a $\bK$-linear space, and let $f:X\times Y\to\bK$ be a a $\bK$-bilinear map.
We call $\bK$-linear space $\{x\in X|\forall y\in Y,f(x,y)=0\}$
(resp. $\{y\in Y|\forall x\in X,f(x,y)=0\}$) of $X$ (resp. $Y$) 
{\it{the left kernel}} (resp. {\it{the right kernel}}).
We call $f$
{\it{non-degenerate}} if the right and left kernels are zero-dimensional.

\begin{lemma} \label{lemma:essngeg}
Let $\BYcnUopflat:=\oplus_{\lambda\in\BYcnfkAoriginal}\bK\BYcntrK_\lambda
(=\BYcnUo\cap\BYcnUpflat)$
and $\BYcnUomflat:=\oplus_{\lambda\in\BYcnfkAoriginal}\bK\BYcntrL_\lambda(=\BYcnUo\cap\BYcnUmflat)$.
Then $\BYcnvartheta_{|\BYcnUopflat\times\BYcnUomflat}$ is non-degenerate.
In particular, $\BYcnvartheta$ is non-degenerate.
\end{lemma}
{\it{Proof.}} Let $\BYcnvartheta^\prime:=\BYcnvartheta_{|\BYcnUopflat\times\BYcnUomflat}$.
Let $\BYcnfkAoriginal^\prime:=\oplus_{i=1}^\BYcnl\bZ\BYcnal_i$
and $\BYcnfkAoriginal^{\prime\prime}:=\oplus_{j=1}^\BYcnl\bZ\BYcnal_{j+\BYcnl}$.
So $\BYcnfkAoriginal=\BYcnfkAoriginal^\prime\oplus\BYcnfkAoriginal^{\prime\prime}$.
Let 
$(\BYcnUomflat)^\prime:=
\oplus_{\lambda\in\BYcnfkAoriginal^\prime}\bK\BYcntrL_\lambda$
and
$(\BYcnUopflat)^{\prime\prime}:=
\oplus_{\mu\in\BYcnfkAoriginal^{\prime\prime}}\bK\BYcntrK_\mu$. 
Let $R$ (resp. $R^\prime$) be the left (resp. right) kernel of $\BYcnvartheta^\prime$.
Let $r\in R$. Write  $r=\sum_{t=1}^zr_t$
for some
$z\in\bN$ and $\lambda_t\in\BYcnfkAoriginal^\prime$
($t\in\BYcnfkJ_{1,z}$) with 
$r_t\in\BYcntrK_{\lambda_t}(\BYcnUopflat)^{\prime\prime}$ ($t\in\BYcnfkJ_{1,z}$).
For $t\in\BYcnfkJ_{1,z}$,
write $\lambda_t=\sum_{i=1}^\BYcnl x^{(t)}_i\BYcnpial_i$ with $x^{(t)}_i\in\bZ$.
Since $\BYcnvartheta^\prime(r,h\BYcntrL_{\BYcnal_{\BYcnl+1}}^{y_1}\cdots\BYcntrL_{\BYcnal_{2\BYcnl}}^{y_\BYcnl})=0$
and
$\BYcnvartheta^\prime(r_t,h\BYcntrL_{\BYcnal_{\BYcnl+1}}^{y_1}\cdots\BYcntrL_{\BYcnal_{2\BYcnl}}^{y_\BYcnl})=\BYcnbrvq^{x^{(t)}_1y_1+\cdots+x^{(t)}_\BYcnl y_\BYcnl}
\BYcnvartheta^\prime(r_t,h)$ 
for all $h\in(\BYcnUomflat)^\prime$ and all $y_i\in\bZ$ ($i\in\BYcnfkJ_{1,\BYcnl}$),
we have
$r_t\in R$ for all $t\in\BYcnfkJ_{1,z}$.
Since $\BYcnvartheta^\prime(r_t,(\BYcnUomflat)^\prime)=\{0\}$,
we have $r_t=0$. Hence $r=0$.
Hence $R=\{0\}$. Similarly we have $R^\prime=\{0\}$.
\hfill $\Box$
\begin{equation*}
\begin{array}{l}
\mbox{We have the Hopf algebra isomorphism} \\ 
\DPD(\BYcnvartheta,\BYcnUp,(\BYcnUm)^{\mathrm{op}})\to\BYcnU,\,\,X\otimes Y\to XY.
\end{array}
\end{equation*}

From now on until the end of Section~\ref{section:CSCent}, we identify
$\BYcnU$ with $\DPD$, let the actions $\ract$, $\lact$, $\duact$ of $\BYcnU$
be the ones defined for $\BYcnU$ by identifying $\BYcnvartheta$ with $\Dpli$.

\subsection{Central elements from finite dimensional modules 
under the assumption~\eqref{eqn:asstilq}}
\label{subsection:cstSkewce}
Keep the notation as in Subsection~\ref{subsection:ExtUtbhmtPi}.
We have assumed \eqref{eqn:asstilq}.
From now on until the end of this subsection
\begin{equation*}
\mbox{fix $\lambda$, $\mu\in\BYcnfkAoriginal$,}
\end{equation*}
and fix a non-zero $\BYcnU$-module $\brvVlm$
such that there exists $\brvvlm\in\brvVlm$
with 
\begin{equation*}
\mbox{$\trK_{\tlambda^\prime}\trL_{\tmu^\prime}\brvvlm:=
\BYcnchi(\lambda,\mu^\prime)\BYcnchi(\tlambda^\prime,\tmu)\brvvlm$
($\lambda^\prime$, $\mu^\prime\in\BYcnfkAoriginal$)
and $\trF_i\brvvlm:=0$ ($i\in\fkI$).}
\end{equation*} 
Define the $\bK$-algebra homomorphism
$\brva:\BYcnU\to\rmEnd_\bK(\brvVlm)$
by $\brva(X)(h):=Xh$ ($X\in\BYcnU$, $h\in\brvVlm$), i.e.,
$(\brva,\brvVlm)$ is the $\bK$-algebra representation of $\BYcnU$ associated with $\brvVlm$.
By Lemma~\ref{lemma:essngeg}, we have
$\brvVlm=\oplus_{\nu\in\BYcnfkAoriginalpip}\BYcnUm_{-\nu}\brvvlm$, as a $\bK$-linear space,
and 
\begin{equation*}
\BYcnUm_{-\nu}\brvvlm=\{\,v\in\brvVlm\,|\,\trK_{\tlambda^\prime}\trL_{\tmu^\prime}v
=\BYcnchi(\tlambda+\nu,\tmu^\prime)\BYcnchi(\tlambda^\prime,\tmu-\nu)v\,\,
(\tlambda^\prime, \tmu^\prime\in\BYcnfkAoriginal)\,\}.
\end{equation*}
For $\nu\in\BYcnfkAoriginalpip$, let $\brvm_\nu:=\dim\trUm_{-\nu}\brvvlm$.
\begin{equation*}
\mbox{Assume $\dim_\bK\brvVlm<\infty$.}
\end{equation*}
Then $\dim_\bK\brvVlm=\sum_{\nu\in\BYcnfkAoriginalpip}\brvm_\nu$.
For $k\in\rmEnd_\bK(\brvVlm)$, let $\rmTrbrvVlm(k)\in\bK$ denote the trace of $k$.
Define $\brvrho\in\rmEnd_\bK(\brvVlm)$ by 
$\brvrho(Y\brvvlm):=\BYcnhrhochipi(\nu)Y\brvvlm$
for $\nu\in\BYcnfkAoriginalpip$ and $Y\in\BYcnUm_{-\nu}$.
Since
$\BYcnS^2(X)=\BYcnhrhochipi(-\lambda^\prime)X$ ($\lambda^\prime\in\BYcnfkAoriginal$,
$X\in\BYcnU_{\lambda^\prime}$),
we have
\begin{equation*}
\trS^2(Y_1)\brvrho(Y_2\brvvlm)=\brvrho(Y_1Y_2\brvvlm)\quad(Y_1,\,Y_2\in\BYcnU).
\end{equation*}
Define $\brvf\in\BYcnUdual$ by $\brvf(X):=\rmTrbrvVlm(\brva(X)\circ\brvrho)$
($X\in\BYcnU$).
\begin{lemma}\label{eqn:lmmefX} We have
\begin{equation}\label{eqn:efX}
\brvf\duact X=\trvarepsilon(X)\brvf\quad(X\in\BYcnU).
\end{equation}
\end{lemma}
{\it{Proof.}} Let $X$, $Y\in\BYcnU$.
Write $\trDelta(X)=\sum_r X^{(1)}_r\otimes X^{(2)}_r$.
Then
\begin{equation*}
\begin{array}{lcl}
(\brvf\duact X)(Y) &=& \brvf(X\lact Y) \\
&=& \sum_r\brvf(X^{(1)}_r Y\trS(X^{(2)}_r)) \\
&=& \sum_r\rmTrbrvVlm(\brva(X^{(1)}_r Y\trS(X^{(2)}_r))\circ\brvrho) \\
&=& \sum_r\rmTrbrvVlm(\brva(Y\trS(X^{(2)}_r))\circ\brvrho\circ\brva(X^{(1)}_r)) \\
&=& \sum_r\rmTrbrvVlm(\brva(Y\trS(X^{(2)}_r))\circ\brva(\trS^2(X^{(1)}_r))\circ\brvrho) \\
&=& \sum_r\rmTrbrvVlm(\brva(Y\trS(X^{(2)}_r)\trS^2(X^{(1)}_r))\circ\brvrho) \\
&=& \sum_r\rmTrbrvVlm(\brva(Y\trS(\trS(X^{(1)}_r)X^{(2)}_r))\circ\brvrho) \\
&=& \trvarepsilon(X)\rmTrbrvVlm(\brva(Y)\circ\brvrho) \\
&=& \trvarepsilon(X)\brvf(Y).
\end{array}
\end{equation*} This completes the proof.
\hfill $\Box$
\newline\par
Define the $\bZ$-module homomorphism 
$\BYcnomegabrzero:\BYcnfkAoriginalpi\to\bKt$ 
by $\BYcnomegabrzero(\BYcnfkAoriginalpi):=\{1\}$.
Then we can easily see the equations below.
\begin{equation}\label{eqn:twoeqs}
\begin{array}{lcl}
\BYcnprtrZomchipibrzero 
& = & \{\,X\in\BYcnU\,|\,XY=YX\,(Y\in\BYcnU)\,\} \\.
& = & \{\,X\in\BYcnU\,|\,\trvarepsilon(Y) X=X\ract Y\,(Y\in\BYcnU)\,\}. 
\end{array}
\end{equation}

\begin{proposition}\label{proposition:prkgfla} We have
\begin{equation} \label{eqn:eqkgflad}
\brvf\in\BYcnDPDBYcnU(\BYcnprtrZomchipibrzero)
\end{equation} and
\begin{equation}\label{eqn:eqkgfla}
\BYcnHCbhmpiomegabrzero(\BYcnDPDBYcnU^{-1}(\brvf))=\sum_{\nu\in\BYcnfkAoriginalpip}
\BYcnhrhochipi(\nu)\brvm_\nu\trK_{\tlambda+\nu}\trL_{\tmu-\nu}.
\end{equation}
\end{proposition}
{\it{Proof.}}
For $\nu\in\BYcnfkAoriginalpip$,
let $m_\nu:=\dim\BYcnUp_\nu$,
let  $\{\,X_{\nu,x}\,|\,x\in\fkJ_{1,m_\nu}\,\}$
and $\{\,Y_{-\nu,y}\,|\,y\in\fkJ_{1,m_\nu}\,\}$
be $\bK$-base of $\trUp_\nu$
and $\BYcnUm_{-\nu}$ respectively
such that $\brvtrvt(X_{\nu,x},Y_{-\lambda^\prime,y})=\krdel_{\nu,\lambda^\prime}\krdel_{x,y}$
($\lambda^\prime\in\BYcnfkAoriginalpip$), and let
$\brvY_{-\nu,z}\in\trUm_{-\nu}$ ($z\in\fkJ_{1,\brvm_\nu}$)
be such that $\{\,\brvY_{-\nu,z}\brvvlm\,|\,z\in\fkJ_{1,\brvm_\nu}\,\}$ is a $\bK$-basis of $\trUm_{-\nu}\brvvlm$.
For $\lambda^\prime\in\BYcnfkAoriginalpip$ and $z\in\fkJ_{1,\brvm_{\lambda^\prime}}$. 
Define $\brvt_{-\lambda^\prime,z}\in\brvVlmdual$
by $\brvt_{-\lambda^\prime,z}(\brvY_{-\nu,y}\brvvlm)=\krdel_{\lambda^\prime,\nu}\krdel_{z,y}$.

For $\lambda^\prime$, $\mu^\prime\in\BYcnfkAoriginalpip$,
$x\in\fkJ_{1,m_{\lambda^\prime}}$, $y\in\fkJ_{1,m_{\mu^\prime}}$ and 
$\lambda^{\prime\prime}$, $\mu^{\prime\prime}\in\BYcnfkAoriginal$, we have
\begin{equation}\label{eqn:impcal}
\begin{array}{l}
\brvf(X_{\lambda^\prime,x}\trK_{\lambda^{\prime\prime}}Y_{-\mu^\prime,y}\trL_{\mu^{\prime\prime}}) \\
\quad = \sum_{\nu\in\BYcnfkAoriginalpip}\sum_{z=1}^{m_\nu}
\brvt_{-\nu,z}(X_{\lambda^\prime,x}\trK_{\lambda^{\prime\prime}}Y_{-\mu^\prime,y}\trL_{\mu^{\prime\prime}}
\brvrho(\brvY_{-\nu,z}\brvvlm)) \\
\quad = \sum_{\nu\in\BYcnfkAoriginalpip}\sum_{z=1}^{m_\nu}
\BYcnhrhochipi(\nu)
\brvt_{-\nu,z}(X_{\lambda^\prime,x}\trK_{\lambda^{\prime\prime}}Y_{-\mu^\prime,y}\trL_{\mu^{\prime\prime}}
\brvY_{-\nu,z}\brvvlm) \\
\quad = \sum_{\nu\in\BYcnfkAoriginalpip}\sum_{z=1}^{m_\nu}
\BYcnhrhochipi(\nu)
\BYcnchi(\lambda^{\prime\prime},-\mu^\prime-\nu+\mu)
\BYcnchi(\nu+\lambda,\mu^{\prime\prime}) \\
\quad\quad\quad\quad\quad\quad\quad\quad
\cdot\brvt_{-\nu,z}(X_{\lambda^\prime,x}Y_{-\mu^\prime,y}
\brvY_{-\nu,z}\brvvlm) \\
\quad = \krdel_{\lambda^\prime,\mu^\prime}\sum_{\nu\in\BYcnfkAoriginalpip}\sum_{z=1}^{m_\nu}
\BYcnhrhochipi(\nu)
\BYcnchi(\lambda^{\prime\prime},-\lambda^\prime-\nu+\mu)
\BYcnchi(\nu+\lambda,\mu^{\prime\prime}) \\
\quad\quad\quad\quad\quad\quad\quad\quad
\cdot\brvt_{-\nu,z}(X_{\lambda^\prime,x}Y_{-\lambda^\prime,y}
\brvY_{-\nu,z}\brvvlm) \\
\quad = \krdel_{\lambda^\prime,\mu^\prime}\sum_{\nu\in\BYcnfkAoriginalpip}\sum_{z=1}^{m_\nu}
\BYcnhrhochipi(\nu)
\brvt_{-\nu,z}(X_{\lambda^\prime,x}Y_{-\lambda^\prime,y}
\brvY_{-\nu,z}\brvvlm) \\
\quad\quad\quad\quad\quad\quad\quad\quad
\cdot\BYcnDPDBYcnU(Y_{-\lambda^\prime,x}\trL_{-\lambda^\prime-\nu+\mu}X_{\mu^\prime,y}
\trK_{\nu+\lambda}) 
(X_{\lambda^\prime,x}\trK_{\lambda^{\prime\prime}}Y_{-\mu^\prime,y}\trL_{\mu^{\prime\prime}}) \\
\quad = \sum_{\nu\in\BYcnfkAoriginalpip}\sum_{z=1}^{m_\nu}
\BYcnhrhochipi(\nu)
\brvt_{-\nu,z}(X_{\lambda^\prime,x}Y_{-\lambda^\prime,y}
\brvY_{-\nu,z}\brvvlm) \\
\quad\quad\quad\quad\quad\quad\quad\quad
\cdot\BYcnDPDBYcnU(Y_{-\lambda^\prime,x}\trL_{-\lambda^\prime-\nu+\mu}X_{\lambda^\prime,y}
\trK_{\nu+\lambda}) 
(X_{\lambda^\prime,x}\trK_{\lambda^{\prime\prime}}Y_{-\mu^\prime,y}\trL_{\mu^{\prime\prime}})  \\
\quad = \sum_{\nu,\nu^\prime\in\BYcnfkAoriginalpip}
\sum_{x^\prime,y^\prime=1}^{m_{\nu^\prime}}
\sum_{z=1}^{m_\nu}
\BYcnhrhochipi(\nu)
\brvt_{-\nu,z}(X_{\nu^\prime,x}Y_{-\nu^\prime,y}
\brvY_{-\nu,z}\brvvlm) \\
\quad\quad\quad\quad\quad\quad\quad\quad
\cdot\BYcnDPDBYcnU(Y_{-\nu^\prime,x^\prime}\trL_{-\nu^\prime-\nu+\mu}X_{\nu^\prime,y^\prime}
\trK_{\nu+\lambda}) 
(X_{\lambda^\prime,x}\trK_{\lambda^{\prime\prime}}Y_{-\mu^\prime,y}\trL_{\mu^{\prime\prime}}) \\
\quad = \BYcnDPDBYcnU(\sum_{\nu,\nu^\prime\in\BYcnfkAoriginalpip}
\sum_{x^\prime,y^\prime=1}^{m_{\nu^\prime}}
\sum_{z=1}^{m_\nu}
\BYcnhrhochipi(\nu)
\brvt_{-\nu,z}(X_{\nu^\prime,x}Y_{-\nu^\prime,y}
\brvY_{-\nu,z}\brvvlm) \\
\quad\quad\quad\quad\quad\quad\quad\quad
\cdot Y_{-\nu^\prime,x^\prime}\trL_{-\nu^\prime-\nu+\mu}X_{\nu^\prime,y^\prime}
\trK_{\nu+\lambda}) 
(X_{\lambda^\prime,x}\trK_{\lambda^{\prime\prime}}Y_{-\mu^\prime,y}\trL_{\mu^{\prime\prime}}).
\end{array}
\end{equation}
Hence $\brvf\in\rmIm\BYcnDPDBYcnU$ and 
\begin{equation}\label{eqn:impcaldash}
\begin{array}{lcl}
\BYcnDPDBYcnU^{-1}(\brvf) & = &
\sum_{\nu,\nu^\prime\in\BYcnfkAoriginalpip}
\sum_{x^\prime,y^\prime=1}^{m_{\nu^\prime}}
\sum_{z=1}^{m_\nu}
\BYcnhrhochipi(\nu)
\brvt_{-\nu,z}(X_{\nu^\prime,x}Y_{-\nu^\prime,y}
\brvY_{-\nu,z}\brvvlm) \\
&  & \quad\quad\quad\quad 
\cdot Y_{-\nu^\prime,x^\prime}\trL_{-\nu^\prime-\nu+\mu}X_{\nu^\prime,y^\prime}
\trK_{\nu+\lambda}.
\end{array}
\end{equation}
By \eqref{eqn:interOmequ}, \eqref{eqn:efX} and \eqref{eqn:twoeqs},
we have $\BYcnDPDBYcnU^{-1}(\brvf)\in\BYcnprtrZomchipibrzero$,
which implies \eqref{eqn:eqkgflad}.
By \eqref{eqn:impcaldash}, we have \eqref{eqn:eqkgfla}.
This completes the proof.
\hfill $\Box$
\newline\par
Let $\BYcnUopi=\BYcnUopichipi:=\oplus_{\nu,\nu^\prime\in\BYcnfkAoriginalpi}
\bK\trK_\nu\trL_{\nu^\prime}$
and $\BYcnUpi=\BYcnUpichipi:=\rmSpan_\bK(\BYcnUm\BYcnUopi\BYcnUp)$.
Then $\BYcnUpi$ (resp. $\BYcnUopi$) is a $\bK$-subalgebra of 
$\BYcnU$ (resp. $\BYcnUo$).
\begin{lemma}\label{lemma:BYom}
Assume $\lambda$, $\mu-\al_{\BYcnl+1}\in\BYcnfkAoriginalpi$.
Let $Z:=\BYcnDPDBYcnU^{-1}(\brvf)\trL_{-\al_{\BYcnl+1}}$.
Then $Z\in\BYcnUpi\cap\BYcnprtrZomchipi$. 
\end{lemma}
{\it{Proof.}}
This easily follows from \eqref{eqn:impcaldash}.
\hfill $\Box$

\subsection{Conjecture on a basis of $\BYcnprtrZomchipi$}
\label{subsection:conjecture}
In this subsection, assume $\BYcnfkAoriginal=\BYcnfkAoriginalpi$.
Then $\BYcnU=\BYcnUpichipi$.
Let $\BYcnomega:\BYcnfkAoriginal\to\bKt$ be a $\bZ$-module homomorphism.
For a $\bK$-algebra homomorphism
$\Lambda:\BYcnUo\to\bK$, there exists a simple
$\BYcnU$-module $\mclL(\Lambda)$
such that there exists $v_\Lambda\in\mclL(\Lambda)\setminus\{0\}$
such that $\mclL(\Lambda)=\oplus_{\nu\in\BYcnfkAoriginalpip}\BYcnUm_{-\nu}v_\Lambda$,
$Zv_\Lambda=\Lambda(Z)v_\Lambda$
and $\trE_iv_\Lambda=0$ ($i\in\BYcnfkI$).

For $\tlambda$, $\tmu\in\tfkAPi$, define the $\bK$-algebra homomorphism
$\BYcnLamlmnewpr:\BYcnUo\to\bK$
by
\begin{equation}\label{eqn:mainweight}
\BYcnLamlmnewpr(\trK_{\lambda^\prime}\trL_{\mu^\prime}):
=\BYcnchi(\lambda,\mu^\prime)
\BYcnchi(\lambda^\prime,\mu)\BYcnomega(\lambda^\prime)
\quad(\lambda^\prime, \mu^\prime\in\BYcnfkAoriginal).
\end{equation} Let
\begin{equation}\label{eqn:mainweightd}
\BYcnFinnewpr:=\{(\lambda,\tmu)\in\BYcnfkAoriginal\times\BYcnfkAoriginal|\dim\mclL(\BYcnLamlmnewpr)<\infty\}.
\end{equation}
\begin{theorem}\label{theorem:mainth}
Let $(\lambda,\mu)\in\BYcnFinnewpr$.
Then
there exists a $\BYcnFinnewprZlm\in\BYcnprtrZomchipi$
such that $\BYcnHCbhmpiomega(\BYcnFinnewprZlm)=
\sum_{\nu\in\BYcnfkAoriginalpip}\BYcnhrhochipi(\nu)
\brvm_\nu\trK_{\lambda+\nu}\trL_{\mu-\nu}$,
where $\brvm_\nu:=\dim\BYcnUm_{-\nu}v_{\BYcnLamlmnewpr}$.
\end{theorem}
{\it{Proof.}} The claim follows from 
Proposition~\ref{proposition:prkgfla} and 
Lemmas~\ref{lemma:BYcniclmm} and \ref{lemma:BYom}.
\hfill $\Box$ 
\newline\par
Now we state the conjecture below.
\newline\par
\begin{equation}\label{eqn:Conjecture}
\mbox{{\bf{Conjecture.}} {\it{The set $\{\,\BYcnFinnewprZlm\,|\,(\lambda,\mu)\in\BYcnFinnewpr\,\}$
must be a $\bK$-basis of $\BYcnprtrZomchipi$.}}}
\end{equation}
\newline\par
\begin{remark}\label{remark:SV2011}
The above conjecture fits into \cite{SV2011}.
\end{remark}

\subsection{$\bZ/3\bZ$-quantum group}
Assume that the characteristic of $\bK$ is not $2$ or $3$.
Assume $\BYcnl=2$. Then $\BYcnfkI=\BYcnfkJ_{1,2}$.
Assume $\BYcnfkAoriginal=\BYcnfkAoriginalpi$.
Let $\sangaV$ be a two dimensional $\bR$-linear space such that
$\BYcnfkAoriginal\subset\sangaV$
and $\{\BYcnpial_1,\BYcnpial_2\}$ is an $\bR$-basis of $\sangaV$.
Namely
$\BYcnfkAoriginal=\BYcnfkAoriginalpi=\bZ\BYcnpial_1\oplus\bZ\BYcnpial_2\subset\sangaV
=\bR\BYcnpial_1\oplus\bR\BYcnpial_2$.
Let $\sangazeta\in \bK$ be such that $\sangazeta^2+\sangazeta+1=0$.
Let $\sangaq\in \bK\setminus\{0,1,\sangazeta\}$.
\begin{equation}\label{eqn:Zthreechi}
\mbox{Let $\BYcnchi(\al_1,\al_1)=\sangazeta$,
$\BYcnchi(\al_2,\al_2)=\sangaq$
and $\BYcnchi(\al_1,\al_2)\BYcnchi(\al_2,\al_1)=\sangaq^{-1}$.}
\end{equation}
Let
$\beta_1:=\BYcnpial_1$, $\beta_2:=2\BYcnpial_1+\BYcnpial_2$, $\beta_3:=\BYcnpial_1+\BYcnpial_2$, 
and $\beta_4:=\BYcnpial_2$.
Then $\BYcnRchipip=\{\beta_t|t\in\BYcnfkJ_{1,4}\}$.
Let 
$\gamma_1:=-\beta_4$, 
$\gamma_2:=-\beta_3$, 
$\gamma_3:=-\beta_2$, 
$\gamma_4:=-\beta_1$, 
$\gamma_5=\beta_4$, 
$\gamma_6=\beta_3$, 
$\gamma_7=\beta_2$, and
$\gamma_8=\beta_1$. 
For $t\in\BYcnfkJ_{1,8}$,
let $\gamma^\prime_t:=\gamma_s$
with $s-t-2\in 8\bZ$.
For $h=(h_1,\ldots,h_8)\in\bRgeqo^8$ and $t\in\BYcnfkJ_{1,8}$,
let $\gamma^{(h)}_t:=\sum_{s=1}^t h_t\gamma_t$.
Let $\gamma^{(h)}_0:=0\in \sangaV$.
Let $B_h:=\cup_{t=1}^8 \{\gamma^{(h)}_{t-1}+u\gamma_t|0\leq u\leq h_t\}$,
and $C_h:=\cap_{t=1}^8(\gamma^{(h)}_{t-1}+\bR\gamma_t+\bRlo\gamma^\prime_t)$.

Let $\Lambda:\BYcnUochipi\to\bK$
be a $\bK$-algebra homomorphism.
Let $l_i:=\Lambda(\BYcntrK_{\BYcnpial_i}\BYcntrL_{-\BYcnpial_i})$
($i\in\BYcnfkI)$.
Assume $\dim\mclL(\Lambda)<\infty$.
Let $\brvm_\nu:=\dim\BYcnUm_{-\nu}v_\Lambda$.
($\nu\in\BYcnfkAoriginalpip$).
By \cite[Lemma~6.6]{AYY12} and its proof, we can see
that there exists a unique $h^\Lambda=(h^\Lambda_1,\ldots,h^\Lambda_8)\in\bZgeqo^8$
such that $\gamma^{(h^\Lambda)}_8=0$,
$\{\nu\in\BYcnfkAoriginalpip|\brvm_\nu>0\}\subset (B_{h^\Lambda}\cup C_{h^\Lambda})\cap\BYcnfkAoriginal$,
and $B_{h^\Lambda}\cap\BYcnfkAoriginalpip\subset\{\nu\in\BYcnfkAoriginalpip|\brvm_\nu=1\}$.
By \cite[(4.15)]{AYY12}, we have
$h^\Lambda_2, h^\Lambda_4, h^\Lambda_6, h^\Lambda_8\in\BYcnfkJ_{0,2}$
and, if $\BYcnkpch(\sangaq)\ne 0$, $h^\Lambda_1, h^\Lambda_3, h^\Lambda_5, h^\Lambda_7\in\BYcnfkJ_{0,\BYcnkpch(\sangaq)-1}$.
Assume $\BYcnkpch(\sangaq)=0$.
Then $\mclL(\Lambda)=\{0\}$ or
there exists $(m,n)\in\bZgeqo\times(2+\bZgeqo)$
such that $l_2=\sangaq^m$, $l_1^2l_2=(\sangazeta\sangaq^{-1})^n$,
$h^\Lambda_{2t}=2$ ($t\in\BYcnfkJ_{1,4}$), $h^\Lambda_1=h^\Lambda_5=m$ and $h^\Lambda_3=h^\Lambda_7=n$.
(See also \cite[Theorem~4.1]{Y07}, \cite[Theorem~7.8.~(cK-4)]{AYY12}.)
Let $H_1:=\{\sangaq^s|s\in\bZ\}$ and $H_2:=\{(\sangazeta\sangaq^{-1})^t|t\in\bZ\}$.
Assume $\BYcnkpch(\sangaq)\ne 0$.
If $l_2\in H_1$, then $h^\Lambda_8=h^\Lambda_2$ and $h^\Lambda_6=h^\Lambda_4$.
If $l_1^2l_2\in H_2$, then $h^\Lambda_8=h^\Lambda_6$ and 
$h^\Lambda_2=h^\Lambda_4$.
For example, letting $\bK=\bC$,
$\xi:=\exp({\frac {2\pi i} {15}})\,(\in\bC)$
($i$ means the imaginary unit),
$\zeta:=\xi^5$ and $q:=\xi^2$, 
if $l_1=1$ and  $l_2=q$,
then $h^\Lambda_1=1$, $h^\Lambda_2=h^\Lambda_8=0$, 
$h^\Lambda_3=h^\Lambda_7=4$
$h^\Lambda_4=h^\Lambda_6=1$,
and $h^\Lambda_5=0$.

By the above argument and Theorems~\ref{theorem:previousmain}
and \ref{theorem:mainth}, 
we can easily convince ourself that 
the conjecture \eqref{eqn:Conjecture} for almost all $\BYcnchi$ of \eqref{eqn:Zthreechi}
must be true.

\section{Appendix}\label{section:Appendix}

In this section, we mention a background of Theorem~\ref{theorem:KhPBW}.
Let $\fkI=\BYcnfkJ_{1,\BYcnl}$ be that of Subsections~\ref{subsection:Prel} and \ref{subsection:DEFtU}.
Let $\bXI$ be the set of all maps from $\bN$ to $\fkI$.
Let the equivalence relation $\sim$ on the set $\bXI\times\bZgeqo$
be such that $(f,t) \sim (f^\prime,t^\prime)$ if and only if
$t=t^\prime$ and $f(k)=f^\prime(k)$ for all $k\in\fkJ_{1,t}$.
Let $\bXsimI$ be the quotient set  $(\bXI\times\bZgeqo)/\!\!\sim$ of $\bXI\times\bZgeqo$ by
$\sim$.
Let 
$x_{f,t}$ be the element of $\bXsimI$
such that $(f,t)$ is a representative of it.
Regard $\bXsimI$ as the totally ordered set  $(\bXsimI,\leq)$ such that 
$x_{f,t}>x_{f^\prime,t^\prime}$ if $t\geq 1$, $t^\prime\geq 1$, $f(s)=f^\prime(s)$ $(s\in\fkJ_{1,t^{\prime\prime}-1}$ and 
$f(t^{\prime\prime})<f^\prime(t^{\prime\prime})$ for some $t^{\prime\prime}\in\fkJ_{1,t}\cap\fkJ_{1,t^\prime}$
or if  $t<t^\prime$ and $f(s)=f^\prime(s)$ $(s\in\fkJ_{1,t})$. (This is the lexicographical order in the sense of \cite[Subsection~1.2]{Kh15}.)
We call $x_{f,t}$ {\it{a word}}.
We call $x_{f,0}$ {\it{an empty  word}}.
For $(f,t)\in\bXI\times\bZgeqo$, define $f^{+,t}\in\bXI$ by $f^{+,t}(k):=f(t+k)$
$(k\in\bN)$.
We also regard $\bXsimI$ as the monoid such that its unit is $x_{f,0}$
and the multiplication is defined in the way that $x_{f,t}x_{f^{+,t},k}=x_{f,t+k}$.
For $i\in\fkI$, let $x_i:=x_{f,1}$ for which $f$ is such that $f(1)=i$.
Then  $\bXsimI$ can also be regarded as the free monoid generated by $x_i$
$(i\in\fkI)$.

Let $\bXsimI^\prime:=\{x_{f,t}\in\bXsimI|t\geq 1\}$ and $\bXsimI^{\prime\prime}:=\{x_{f,t}\in\bXsimI|t\geq 2\}$.

We say that $u=x_{f,t}\in\bXsimI$ is {\it{a standard word}} (or {\it{a Lyndon-Shishov word}}) if $t\in\fkJ_{0,1}$
or if $u\in\bXsimI^{\prime\prime}$ and $u>wv$ for all $v$, $w\in\bXsimI^\prime$ 
with $vw=u$.
Let $\bLSsimI$ be the set of all standard words. 
Let $\bLSsimI^\prime:=\bLSsimI\cap\bXsimI^\prime$.
and $\bLSsimI^{\prime\prime}:=\bLSsimI\cap\bXsimI^{\prime\prime}$.
It is well-known \cite[Corollary1.1]{Kh15} that
\begin{equation*}
\forall u\in \bLSsimI^{\prime\prime}, \exists v, \exists w\in \bLSsimI^\prime\,
{\mathrm{s.t.}}\, u=vw, v>w.
\end{equation*}
For $u\in \bLSsimI^{\prime\prime}$,
let ${\dot{u}}$, ${\ddot{u}}\in\bLSsimI^\prime$
be such that 
$u={\dot{u}}{\ddot{u}}$
and there exists $y\in\bXsimI$
with $v={\dot{u}}y$ for all $v$, $w\in \bLSsimI^\prime$ satisfying $u=vw$ 
(\cite[Theorem~1.1]{Kh15} implies ${\dot{u}}>{\ddot{u}}$ and $v>w$).
Let $\bK\bXsimI$ be the $\bK$-linear space such that
$\bXsimI$ is its $\bK$-base.
Regard $\bK\bXsimI$ as the associative $\bK$-algebra such that
$\bXsimI$ is a submonoid of $\bK\bXsimI$.
Then  $\bK\bXsimI$ can also be regarded as the free $\bK$-algebra generated by $x_i$
$(i\in\fkI)$.   

Let $\BYcnchi:\BYcnfkAoriginal\times\BYcnfkAoriginal\to\bKt$,
$\al_i$ $(i\in\BYcnfkI)$ and $\BYcnfkAoriginalpip$ be as in Subsections~\ref{subsection:Prel} and \ref{subsection:DEFtU}.
For $x_{f,t}\in\bXsimI$, let $\theta(x_{f,t}):=\sum_{k=1}^t\BYcnal_{f(k)}(\in\BYcnfkAoriginalpip)$.
For $u=x_{f,t}\in \bLSsimI$, define $[u]\in\bK\bXsimI$
by 
\begin{equation*}
[u]:=\left\{\begin{array}{ll}
u & \quad\mbox{if $t\in\fkJ_{0,1}$}, \\
\mbox{$[{\dot{u}}][{\ddot{u}}]$}
-\chi(\theta({\ddot{u}}),\theta({\dot{u}}))^{-1}
[{\ddot{u}}][{\dot{u}}]
& \quad\mbox{if $t\geq 2$}.
\end{array}\right.
\end{equation*}
We call $[u]$ a {\it{ super-letter}}.

\begin{theorem}\label{theorem:KhaA}
{\rm{(}}See \cite[Theorems~1.1, 2.6 and Lemma~2.6]{Kh15}{\rm{)}}
{\rm{(1)}} The set 
\begin{equation}\label{eqn:defbB}
\{1\}\cup\{[u_1][u_2]\cdots [u_k]
|k\in\bN,\,u_r\in\bLSsimI^\prime\,
(r\in\fkJ_{1,k}),
u_r\leq u_{r+1}\,(r\in\fkJ_{1,k-1})\}
\end{equation} is a $\bK$-basis of $\bK\bXsimI$. 
Moreover, for $v_1$, $v_2\in\bLSsimI^\prime$ with $v_1>v_2$,
$[v_1][v_2]$ is a linear combination of elements $[u_1]^{n_1}[u_2]^{n_2}\cdots[u_k]^{n_k}$
in \eqref{eqn:defbB} with
$v_2\leq u_1$, $v_1\leq u_k$, and $\theta(u_1^{n_1}u_2^{n_2}\cdots u_k^{n_k})
=\theta(v_1)+\theta(v_2)$.
\newline
{\rm{(2)}} The same claim as that of {\rm{(1)}} with the monomials $u_1u_2\cdots u_k$ 
in place of $[u_1][u_2]\cdots [u_k]$ is true.
\end{theorem}

Let $\acuteBYcnUpflat$ be a Hopf $\bK$-algebra
generated by $\acuteBYcntrK_\lambda$
$(\lambda\in\BYcnfkAoriginal)$ and $\acuteBYcntrE_i$
$(i\in\BYcnfkI)$
satisfying the same equations as those of 
\eqref{eqn:relone} and \eqref{eqn:defHopf} with 
$\acuteBYcntrK_\lambda$ and  $\acuteBYcntrE_i$
in place of
$\BYcntrK_\lambda$ and  $\BYcntrE_i$
respectively.
Let $\acuteBYcnUpflato$ (resp. $\acuteBYcnUp$)
be the $\bK$-subalgebra of $\acuteBYcnUpflat$
generated by $\acuteBYcntrK_\lambda$
$(\lambda\in\BYcnfkAoriginal)$
(resp. $\acuteBYcntrE_i$
$(i\in\BYcnfkI)$).
Note $\acuteBYcnUpflat=\rmSpan_\bK(\acuteBYcnUp\acuteBYcnUpflato)$.
Let $\acuteG:=\{\acuteBYcntrK_\lambda|\lambda\in\BYcnfkAoriginal\}$.
Assume that $\acuteG$
is a $\bK$-basis of $\acuteBYcnUpflato$, and that
$\acuteBYcntrK_\lambda\ne\acuteBYcntrK_\mu$ for $\lambda\ne\mu$.
Assume that $\acuteG=\{g\in\acuteBYcnUpflat|\BYcnDelta(g)=g\otimes g\}$. 
(Such Hopf algebra is called a {\it{Character Hopf algebra}}, 
see \cite[Defnnition~1.11]{Kh15}.)
Let $\Gamma=(\Gamma,\preceq)$ be a well-ordered additive (commutative) monoid,
see \cite{Kob04} for this terminology.
\begin{example}
{\rm{
Let $n\in\bN$. Define the well-order $\preceq$ on $\bZgeqo^n$ in the way that
for $k=(k_1,\ldots, k_n)$, $r=(r_1,\ldots, r_n)\in\bZgeqo^n$, 
we have $k\prec r$ if and only if $\sum_{t=1}^nk_t<\sum_{t=1}^nr_t$
or there exists $b\in\BYcnfkJ_{1,n}$ with
$k_b<r_b$ and $k_{b^\prime}=r_{b^\prime}$ for all $b^\prime\in\BYcnfkJ_{1,b-1}$.
Then $(\bZgeqo^n,\preceq)$ is a well-ordered additive monoid.
Note that $\BYcnfkAoriginalpip$ is a well-ordered additive monoid
in a natural sense.
}}
\end{example}

Let $D:\bXsimI\to\Gamma$ be a monoid homomorphism such that
$0\notin D(\bXsimI^\prime)$.
Also, let $D([u_1][u_2]\cdots [u_k])$
mean $D(u_1u_2\cdots u_k)$.
Define the $\bK$-algebra epimorphism $\Phi:\bK\bXsimI\to\acuteBYcnUp$ by
$\Phi(x_i):=\acuteBYcntrE_i$ $(i\in\BYcnfkI)$.
For the element $X=[u_1][u_2]\cdots [u_k]$ in \eqref{eqn:defbB},
let $M(X):=u_k$.
Let $\bB$ be the set in \eqref{eqn:defbB}.
For $u\in\bLSsimI^\prime$, we call $[u]$ {\it {hard}} 
if $\Phi([u])\notin
\rmSpan_\bK(\{\Phi(X)|X\in\bB,D(X)=D(u), M(X)<u\}\cup\{g\Phi(Y)|g\in\acuteG, Y\in\bB, D(Y)<D(u)\})$,
see \cite[Definition~~2.6 (see also Lemma~2.6)]{Kh15}.
For a hard super-letter $[u]$,
let $H([u]):=\{k\in\bN|
\Phi([u]^k)\in
\rmSpan_\bK(\{\Phi(X)|X\in\bB,D(X)=D(u^k), M(X)<u\}\cup\{g\Phi(Y)|g\in\acuteG, Y\in\bB, D(Y)<D(u^k)\})\}$.
If $H([u])\ne\emptyset$, let $h([u]):=\min H([u])$.
Otherwise, let $h([u]):=\infty$.
Clearly $h([u])\geq 2$.

\begin{theorem}\label{theorem:KhMain}
{\rm{(\cite[Defintion~2.7, Theorems~2.2 and 2.4 (see also Lemma~2.6)]{Kh15})}}
{\rm{(1)}} Let $h:=h([u])$. If $h<\infty$, then 
$\chi(\theta(u),\theta(u))$ is a primitive $t$-th root of unity for some $t\in\bN$,
and $h=tp^r$ for some $r\in\bZgeqo$, where $p:=\mathrm{Char}(\bK)(\in\{0\}\cup(\bN\setminus\{1\})$.
\newline 
{\rm{(2)}} Let $\bH$ be the set of hard super-letters.
Then the elements 
\begin{equation*}
g\Phi([u_1]^{n_1}[u_2]^{n_2}\cdots[u_k]^{n_k})\quad (g\in\acuteG,[u_i]\in\bH,
n_i<h([u_i]),u_1<\cdots<u_k)
\end{equation*} form a $\bK$-base of $\acuteBYcnUpflat$. 
\newline
{\rm{(3)}} The same claim as that of {\rm{(2)}} with the monomials $u_1^{n_1}u_2^{n_2}\cdots u_k^{n_k}$ 
in place of $[u_1]^{n_1}[u_2]^{n_2}\cdots[u_k]^{n_k}$ is true.
\end{theorem}

Assume that we have the $\bK$-linear isomorphism 
$\acuteBYcntrm:\acuteBYcnUp\otimes\acuteBYcnUpflato\to\acuteBYcnUpflat$ defined
by $\acuteBYcntrm(X\otimes Z):=XZ$. By Theorem~\ref{theorem:KhMain}, we can 
easily see that
\begin{equation}\label{eqn:GenEqKhPBW}
\begin{array}{l}
\mbox{$\{\Phi([u_1]^{n_1}[u_2]^{n_2}\cdots[u_k]^{n_k}) 
|[u_i]\in\bH,
n_i<h([u_i]),u_1<\cdots<u_k\}$} \\ 
\mbox{is a $\bK$-basis of $\acuteBYcnUp$.} 
\end{array}
\end{equation} (The claim of $[u_1]^{n_1}[u_2]^{n_2}\cdots[u_k]^{n_k}$ can be replaced by $u_1^{n_1}u_2^{n_2}\cdots u_k^{n_k}$.)
Let $\acutebH:=\bH\cup\{[u]^{\BYcnkpch(\chi(\theta(u),\theta(u)))p^s}|[u]\in\bH, s\in\bZgeqo,0<\BYcnkpch(\chi(\theta(u),\theta(u)))p^s<h([u])\}$, where $p:=\mathrm{Char}(\bK)$.
Define the total order $<$ on $\acutebH$
in the way that $[u]^x>[v]^y$ if and only if $u>v$ or $u=v$ and $x<y$.
Define the map $\acuteh:\acutebH\to(\bN\setminus\{1\})\cup\{\infty\}$ as follows.
If $\BYcnkpch(\chi(\theta(u),\theta(u)))=0$, let $\acuteh([u]):=\infty\,(=h([u]))$.
Otherwise, let $\acuteh([u]^x):=\BYcnkpch(\chi(\theta(u^x),\theta(u^x)))$.
Then the same claim as \eqref{eqn:GenEqKhPBW}
with $\acutebH$ and $\acuteh$ in place of
$\bH$ and $h$ is true.

Assume that we have a direct sum
$\acuteBYcnUp=\oplus_{\lambda\in\BYcnfkAoriginalpip}\acuteBYcnUp_\lambda$
as a $\bK$-linear space
with
$1\in\acuteBYcnUp_0$,
$\acuteBYcntrE_i\in\acuteBYcnUp_{\al_i}$
and $\acuteBYcnUp_\lambda\acuteBYcnUp_\mu\subset\acuteBYcnUp_{\lambda+\mu}$.
Let $\acuteRchip:=\{\theta(u^x)|[u]^x\in\acutebH\}$.
Define the map 
$\acutephichip:\acuteRchip\to\bN$
by $\acutephichip(\beta):=|\{[u]^x\in\acutebH|\theta(u^x)=\beta\}|$.
Then the pair $(\acuteRchip,\acutephichip)$ is independent from choice
of $D$ (and $\Gamma$).
Recall $\BYcnUpflat$ from the sentence below \eqref{eqn:defHopf}.
Identify $\acuteBYcnUpflat$ with $\BYcnUpflat$
in a natural way.
Then $(\acuteRchip,\acutephichip)$ can be identified with $(\BYcnRchipip,\BYcnvphpip)$
of Theorem~\ref{theorem:KhPBW}. Let $[u]^x\in\acutebH$,
$\al:=\theta([u]^x)$,
$t\in\fkJ_{1,\acutephichip(\al)}$ and $d:=(\al,t)$. Then we have $y(d)\in\fkJ_{0,\acuteh([u]^x)-1}$,
where $y$ is the one of Theorem~\ref{theorem:KhPBW}.

\vspace{1cm}
{Punita Batra, 
Harish-Chandra Research Institute,
Chhatnag Road,
Jhunsi,
Allahabad 211 019,
India}, \newline
{E-mail: batra@hri.res.in}
\vspace{0.5cm}

{Hiroyuki Yamane, Department of Mathematics,
Faculty of Science, University of Toyama,
3190 Gofuku, Toyama-shi, Toyama 930-8555, Japan}, \newline
{E-mail: hiroyuki@sci.u-toyama.ac.jp}

\end{document}